\documentclass{amsart}

\usepackage[latin1]{inputenc}

\usepackage{amsmath,amssymb,amscd,latexsym,color}

\usepackage[dvips]{graphicx}

\usepackage[english]{babel}

\def\Act{A}
\def\blue{\text{b}}
\def\yellow{\text{y}}
\def\green{\text{g}}
\def\bluestar{\text{b}^*}
\def\yellowstar{\text{y}^*}
\def\greenstar{\text{g}^*}

\def\shell{\text{Shell}}
\def\signe{\epsilon}

\newcommand{\communique}{\leftrightarrow}
\newcommand{\communiquep}{\stackrel{p}{\leftrightarrow}}
\newcommand{\communiqueq}{\stackrel{q}{\leftrightarrow}}

\newcommand{\N}{\mathbb{Z}_{+}}

\newcommand{\Z}{\mathbb{Z}}

\newcommand{\Zd}{\mathbb{Z}^d}

\newcommand{\Q}{\mathbb{Q}}

\newcommand{\Pcond}{\overline{\mathbb{P}}}

\newcommand{\Rd}{\mathbb{R}^d}

\renewcommand{\P}{\mathbb{P}}

\newcommand{\Ld}{\mathbb{L}^d}

\newcommand{\Ed}{\mathbb{E}^d}

\newcommand{\Ber}{\text{Ber}}

\newcommand{\as}{\text{ a.s.}}
\newcommand{\io}{\text{ i.o.}}

\renewcommand{\epsilon}{\varepsilon}
\renewcommand{\phi}{\varphi}
\renewcommand{\limsup}{\overline{\lim}}
\renewcommand{\liminf}{\underline{\lim}}

\newcommand{\iid}{\text{i.i.d. }}
\newcommand{\ie}{\emph{i.e. }}

\newcommand{\resp}{\emph{resp. }}

\newcommand{\miniop}[3]{%
\renewcommand{\arraystretch}{0.6}
\begin{array}{c}
{\scriptstyle #1}\\
#2\\
{\scriptstyle #3}
\end{array}
\renewcommand{\arraystretch}{1}}

\newcommand{\1}{1\hspace{-1.3mm}1}

\newcommand{\pcfleche}{\overrightarrow{p_c}}

\begin{document}

{
\newtheorem{theorem}{Theorem}[section]
\newtheorem{conjecture}[theorem]{Conjecture}

}
\newtheorem{lemme}[theorem]{Lemma}
\newtheorem{defi}[theorem]{Definition}
\newtheorem{coro}[theorem]{Corollary}
\newtheorem{rem}[theorem]{Remark}
\newtheorem{prop}[theorem]{Proposition}

\title[Competition and Bernoulli Percolation]{Competition between  growths governed by Bernoulli Percolation}
\date{\today}

{
\author{Olivier Garet}
\author{R{\'e}gine Marchand}
\address{Laboratoire de Math{\'e}matiques, Applications et Physique
Math{\'e}matique d'Orl{\'e}ans UMR 6628\\ Universit{\'e} d'Orl{\'e}ans\\ B.P.
6759\\
 45067 Orl{\'e}ans Cedex 2 France}
\email{Olivier.Garet@univ-orleans.fr}

\address{Institut Elie Cartan Nancy (math{\'e}matiques)\\
Universit{\'e} Henri Poincar{\'e} Nancy 1\\
Campus Scientifique, BP 239 \\
54506 Vandoeuvre-l{\`e}s-Nancy  Cedex France}
\email{Regine.Marchand@iecn.u-nancy.fr}
}

\def\motsclefs{Percolation, first-passage percolation, chemical distance, competition, random growth.}

{
}

\begin{abstract}
We study a competition model on $\Zd$ where the two infections are driven by supercritical Bernoulli percolations with distinct parameters $p$ and $q$. We prove that, for any $q$, there exist at most countably many values of $p<\min{\{q,\pcfleche\}}$ such that coexistence can occur.

\end{abstract}

{\maketitle 
}

\section{Introduction}

Consider two infections, say \emph{blue} and \emph{yellow}, which attempt to conquer, in discrete time, the space $\Zd$. 
At time $0$, all sites are empty but two: one is \emph{active blue}, source of the blue infection, the other one is \emph{active yellow}, source of the yellow infection. To evolve from time $t$ to time $t+1$, the process is governed by the following rules. Each infection is only transmitted by \emph{active} sites of its color to \emph{empty} sites.
Each active site tries to infect each of its empty neighbors, and succeeds with probability $p_{\blue}$ or $p_{\yellow}$, according to its color, blue or yellow. In case of success, the non-occupied site becomes an active site with the color of the infection; otherwise, it remains empty. In any case, the active site becomes a \emph{passive} site of the same color, and can not transmit any infection anymore.
Moreover, we make the following assumptions:
\begin{itemize}
\item  the success of each attempt of contamination at a given time does not depend on the past,
\item the successes of simultaneous attempts of contamination are independent.
\end{itemize}
The first point allows a modelization of this competition model by a homogeneous Markov chain while Markov chains satisfying the second point are sometimes called \emph{Probabilistic Cellular Automata (PCA)}.
Note that if the two initial sources are at an \emph{odd} $\|.\|_1$-distance from each other, no empty site will be infected at the same time by the two distinct infections. To extend the definition of the model to more general initial configurations, we will add some extra rules in the next section. 

Thus the two infections compete to invade space: once a site is colored, it keeps its color for ever and cannot be used by the other infection as a transmitter. As in other competition models, it is natural to ask whether coexistence, \ie unbounded growth of the two infections, can occur. We propose 
the following conjecture:
\begin{conjecture}
\label{Introconj}
If $p_{\blue}=p_{\yellow}>p_c$ then coexistence occurs with positive probability, while if $p_{\blue} \neq p_{\yellow}$ and at least one them is strictly smaller than $\overrightarrow{p_c}$ then coexistence cannot occur.
\end{conjecture}
We will soon see that this competition model is closely linked to Bernoulli bond percolation on $\Zd$, where, as usually, $p_c=p_c(d)$ denotes the critical probability and $\overrightarrow{p_c}=\overrightarrow{p_c}(d)$ the critical probability in the oriented setting. This link will easily explain the fact that if $p_{\blue}<p_c$ -- \resp $p_{\yellow}<p_c$ --  then with probability one the blue -- \resp yellow -- infection dies out. Thus the interesting case is when each infection has a parameter larger than $p_c$. In the same manner, it is not difficult to see, using properties of supercritical oriented Bernoulli percolation, that if both  $p_{\blue}$ and $p_{\yellow}$ exceed $\pcfleche$, then  the
 blue and yellow infections can simultaneously grow unboundedly.

The coexistence statement of the conjecture has already been proved in a previous paper of the authors~\cite{GM-coex}. To precise the second part, and before stating the corresponding main result of this paper, we would like to recall the state of the art in competition problems of this type. A very natural way to obtain a competition model is to extend some well understood one-type interacting particle system in such a way that each infection
behaves like the one-type model does in each region where only one of both 
types is present.
Some famous one-type interacting particle have been considered:
the contact process by Neuhauser~\cite{neuhauser}, the Richardson model by
Häggström and Pemantle~\cite{Haggstrom-Pemantle-1,Haggstrom-Pemantle-2}, or by another way by Kordzakhia and Lalley~\cite{KL}, Deijfen's continuous version
of Richardson model~\cite{deijfen} by Deijfen, Haggström and Bagley~\cite{DHB} 
and Deijfen and Häggström~\cite{DH}. Each of these models actually corresponds to a family of stochastically comparable processes indexed by a continuous parameter and it is natural to ask if coexistence is possible when the two infections are governed by the same parameter -- \resp by different parameters. Note that in all these models, the stochastically comparable processes are governed by exponential families. The following dichotomy seems to emerge.

Either the two infections have the same strength, or same speed of propagation. In this case, coexistence occurs with positive probability: 
it has been proved at first for the two-type Richardson model when $d=2$ by Häggström and Pemantle~\cite{Haggstrom-Pemantle-1} and then extended by Garet and Marchand~\cite{GM-coex} for
a wide class of first-passage percolation models, including the percolation model that is studied here. An alternative proof is also given by Hoffman~\cite{hoffman}. Similarly,  Deijfen and Häggström~\cite{DH} proved the possibility of 
coexistence for Deijfen's continuous version
of Richardson model.
The same result is also proved by Kordzakhia and Lalley~\cite{KL} for their own
extension of Richardson model. Nevertheless, their proof is conditioned
by a difficult and reasonable conjecture on the curvature properties of
the asymptotic shape for Richardson model.

Or one infection is stronger -- or faster than the other one. It is then conjectured that coexistence is not possible. The first and famous result in this direction was done by Häggström and Pemantle~\cite{Haggstrom-Pemantle-2}: they proved that for their model,
coexistence is not possible, except perhaps for a denumerable set for the
ratio of the speeds. The result of   Deijfen, Häggström and Bagley~\cite{DHB}
is submitted to the same irritating restriction rule. 

\begin{figure}
\begin{tabular}{ccc}
\includegraphics[scale=0.8]{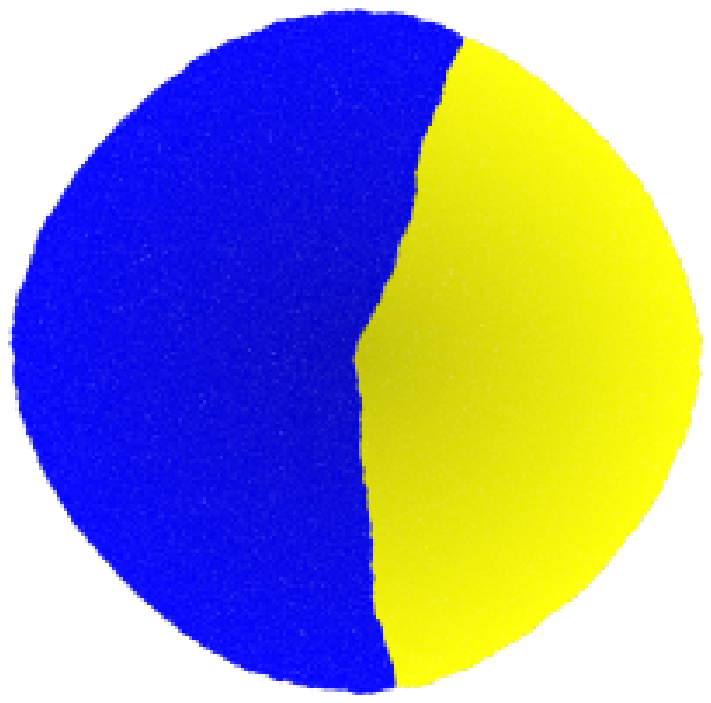} &  &
\includegraphics[scale=0.8]{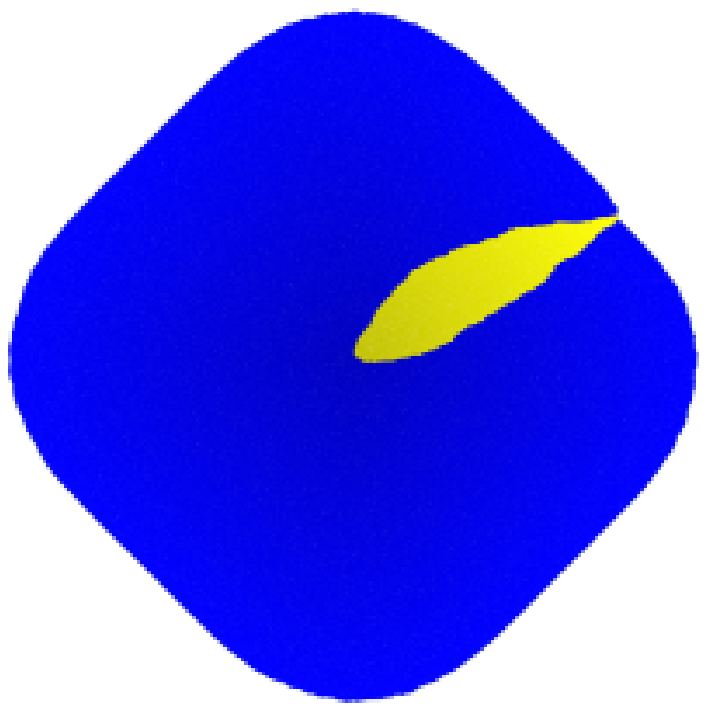}\\
$p_{\blue}=p_{\yellow}=0.6$ & \hspace{-2cm} & $p_{\blue}=0.7\ p_{\yellow}=0.6$
\end{tabular}
\caption{Bernoulli competition in a $4000\times 4000$ grid}
\end{figure}

For our model, we
prove, in this paper, the similar following result:
\begin{theorem}
\label{THEtheorem}
Let $p_{\blue}>p_c$ be fixed:
there exists a denumerable set $\text{Bad}\subset [p_c,\pcfleche)$ such that
for each $p_{\yellow}\in [0,\min\{p_{\blue},\pcfleche\})\backslash \text{Bad}$, 
the probability that both infections infinitely grow is null.
\end{theorem}

Before commenting this result, and to complete the survey, let us mention the recent paper
by Deijfen and Häggström~\cite{DH-nonmon}, where they exhibit graphs where
coexistence occurs for several values for the ratio of the speeds.
This should prevent researchers from unsuccessful attempts to fill the gap 
with the only help of stochastic comparisons. 

In its main lines, the present paper follows the strategy initiated by 
Häggström and Pemantle~\cite{Haggstrom-Pemantle-2}, but it has to overcome some extra difficulties.
Our model also depends on one simple parameter -- the parameter of the related Bernoulli percolation -- which allows coupling and stochastic comparisons. However, note that:
\begin{itemize}
\item The memoryless properties of the exponential laws are lost: one active site tries to infect an empty neighbor only once.
\item Scaling properties of the asymptotic shape in first-passage percolation with exponential times are lost: asymptotic shapes corresponding 
to different values of the parameter are not homothetic anymore.
\end{itemize}

The paper is organized as follows. First, in Section~\ref{Sectiondefinition}, we describe precisely the PCA underlying this competition process, exhibit its reformulation in terms of Bernoulli percolation, and give some related coupling properties and stochastic comparison results. Then,  Section~\ref{Sectionbernoulli} gives a primer of results concerning Bernoulli percolation and the related chemical distance: we particularly recall there the convergence result of the chemical distance in Bernoulli percolation with supercritical parameter $p$ to a norm $\|.\|_p$, and an associated large deviation result.

The first key point of the proof of the main result is the strict comparison of the norms associated to the asymptotic behavior of chemical distance in Bernoulli percolation with different parameters, which will replace the homothetic properties of asymptotic shapes in the case of exponential laws. Section~\ref{Sectionstrictecomp} is devoted to the proof of this result:
\begin{prop} 
\label{Introstrictecomp}
Assume that $p_c<p<\overrightarrow{p_c}$ and $p<q \le 1$. There exists a positive constant $C_{p,q}<1$ such that
$$\forall x \in \Rd \quad \|x\|_q \le C_{p,q}\|x\|_p.$$
\end{prop}
Although the large comparison $\|x\|_q \le \|x\|_p$ is quite natural, the strict comparison will be necessary to ensure, roughly speaking, that in every direction, the stronger infection can take a real advantage and grow strictly faster than the other one.

The second key step is to prove that when coexistence occurs, the global growth of the infected sites is governed by the norm of the weaker infection: denote by $\eta(t)$ is the set of already infected sites at time $t$ and $|A|_p=\sup\{\|x\|_p: \; x\in A\}$. Then  
\begin{prop}
\label{Introslowspeed}
Let $p$ and $q$ be such that $p_c<q\le 1$ and $0\le p<\min\{q,\pcfleche\}$. On the event ``the weak infection survives'', we have almost surely:
$$\miniop{}{\limsup}{t\to +\infty}\frac{|\eta(t)|_{p}}t\le 1.$$
\end{prop} 
The proof of this proposition -- in fact the core of the paper -- is given in Section~\ref{Sectionslowspeed}. It relies both on the previous proposition and on the large deviation result on the set of infected points with respect to the asymptotic shape in the corresponding one infection
model which is recalled in Section~\ref{Sectionbernoulli}.

Finally, in Section~\ref{Section_proof_theo}, we collect all these results to prove the main theorem via coupling results that are in the spirit of Häggström and Pemantle's work~\cite{Haggstrom-Pemantle-2}.

\section{The competition model}
\label{Sectiondefinition}

This section has several goals:
\begin{itemize}
\item to complete the progression rules exposed in the introduction and to define the model for general initial configurations. This will correspond to the artificial introduction of \emph{green} sites.
\item to define the PCA by describing the transition matrix of the homogeneous Markov chain in terms of local rules.
\item to give an alternative description~(\ref{THEdefinition}) in terms of Bernoulli percolation and chemical distance and to prove the equivalence between the two definitions in Lemma~\ref{randomset}. This definition will be the one used in the next sections.
\item Use this last definition to give monotonicity properties in Lemma~\ref{yal} and comparisons properties between the one-type growth and the two-types growth model in Lemma~\ref{comparaison}.
\end{itemize}

Suppose from now on that $p_{\yellow}\le p_{\blue}$, which means that the blue infection is stronger than the yellow one. To complete the description of the model, let us first describe the interface between the two infections via the introduction of \emph{green} sites. A green site is to be understood as a superposition of a blue site and a yellow site. To be coherent with the previous rules, we assume that an active green site transmits to each of its empty neighbors either both infections with probability $p_{\yellow}$, or only the blue infection with probability $p_{\blue}-p_{\yellow}$, or fails in its infection attempts with probability $1-p_{\blue}$; it then becomes a passive green site. 
Note that this rule is quite arbitrary. The necessary part is that a green site transmits to one of its neighbor a yellow -- \resp blue -- infection with probability $p_{\yellow}$ -- \resp $p_{\blue}$ -- and we choose the coupling between these two transmissions to simplify some coupling in the sequel, but it has no real influence on the behavior of the model.

To determine the state at time $t+1$ of an empty site $x$ at time $t$, we then check the types of infections that are transmitted to it: either they are all of the same color, blue or yellow, and $x$ becomes an active site of this color, or they are of both colors, and $x$ becomes an active green site, or no infection is transmitted to $x$, which then remains empty. We can now give the formal definition of the PCA.

\subsection{Definition of the Probabilistic Cellular Automata (PCA)}

\subsubsection*{Definition of the graph $\Ld$}
We endow the set $\Zd$ with the set of edges $\Ed$ between sites of $\Zd$ that are at distance 1 for the Euclidean distance: the obtained graph is denoted $\Ld$. Two sites $x$ and $y$ that are linked by an edge are said to be \emph{neighbors} and this relation is denoted: $x \sim y$.

\subsubsection*{State space}
Let us introduce the set $S=\{0,\blue,\yellow,\green,\bluestar,\yellowstar,\greenstar\}$ of possible states of a site:
$0$ is the state of an empty site, $\blue,\yellow,\green$ -- corresponding respectively to colors blue, yellow and green -- the states
of active sites, and $\bluestar,\yellowstar,\greenstar$ the states of passive colored sites.

In the sequel, we will restrict our Markov chain to start from a configuration with a finite number of colored sites, whence the only configurations appearing during the whole process will also have a finite numbers of colored sites. Our Markov chain will thus live in the following denumerable state set:
$$S^{(\Zd)}=\{\xi\in S^{\Zd}: \; \exists \Lambda\text{ finite},\; \xi_k=0\text{ for }k\in\Zd\backslash\Lambda\}.$$ 

\subsubsection*{Local rules}
To complete the definition of the Markov chain, it only remains to define its transition probabilities, via \emph{local rules}, describing the evolutions of the infections exposed in the introduction.
Define, for $c \in \Act=\{\blue,\yellow,\green\}$, the number $n_x^\xi(c)$ of active neighbors with color $c$ of the site $x\in \Zd$ in the configuration $\xi \in S^{(\Zd)}$: 
$$n_x^\xi(c)=|\{y\in\Zd: \; x \sim y \text{ and }\xi_y=c\}|,$$
and define the probability $p_x^\xi(c, \tilde c)$ that the site $x \in \Zd$, in the configuration $\xi \in S^{(\Zd)}$, swaps from color $c$ to color $\tilde c$:
\begin{itemize}
\item If $x$ is an empty site, \ie if $\xi_x=0$, set:
$$
\left\{
\begin{array}{rcl}
p_x^\xi(0,0)
& = & (1-p_\blue)^{n_x^\xi(\blue)+n_x^\xi(\green)} (1-p_\yellow)^{n_x^\xi(\yellow)} \\
p_x^\xi(0,\yellow)
& = & (1-p_\blue)^{n_x^\xi(\blue)+n_x^\xi(\green)} \left[ 1-(1-p_\yellow)^{n_x^\xi(\yellow)} \right] \\
p_x^\xi(0,\blue)
& = & \left[ 1-(1-p_\blue)^{n_x^\xi(\blue)} \right] (1-p_\yellow)^{n_x^\xi(\yellow)+n_x^\xi(\green)} \\
& & + (1-p_\blue)^{n_x^\xi(\blue)} (1-p_\yellow)^{n_x^\xi(\yellow)} \left[ 1-(1-p_\blue+p_\yellow)^{n_x^\xi(\green)} \right] \\
p_x^\xi(0,\green)
& = & 1- p_x^\xi(0,0)-p_x^\xi(0,\yellow)-p_x^\xi(0,\blue).
\end{array}
\right.
$$
\item If $x$ is an active site, it becomes passive: $\forall c \in \Act, \; p_x^\xi(c, c^*) =  1$.
\item If $x$ is an passive site, it remains passive: $\forall c \in \Act, \; p_x^\xi(c^*, c^*) =  1$.
\item In any other case, the probability is null.
\end{itemize}

\subsubsection*{Transition probabilities}
We can then define the following transition probabilities on the state set $S^{(\Zd)}$:
\begin{equation}
\label{matrix}
\forall (\xi^1,\xi^2)\in S^{(\Zd)}\times S^{(\Zd)}\quad p(\xi^1,\xi^2)=\prod_{x\in\Zd} p_x^{\xi^1}(\xi^1_x,\xi^2_x).
\end{equation}
Note that, as only a finite number of terms differ from~$1$, the previous product is convergent. 

\subsection{Realization of the Markov chain via Bernoulli percolation} 
\label{partbernouill}

The aim of this part is to link this PCA with some natural Bernoulli percolation structures on $\Zd$, and to give an alternative description of the model in terms of random sets and of a specific first-passage percolation model. We begin with some classical notations of Bernoulli percolation on $\Zd$.

\subsubsection*{Bernoulli percolation}
On the set $\Omega=[0,1]^{\Ed}$ endowed with its Borel $\sigma$-algebra, consider the probability measure $\P=\text{Unif}[0,1]^{\otimes\Ed}$. 
For each $p\in [0,1]$ and $\omega\in\Omega$, 
we denote by $\mathcal G_p(\omega)$ the subgraph of $\Ld$ whose bonds $e$ are $p$-\emph{open}, which means that they are are such
that $\omega_e\le p$.

For $A\subset\Zd$ and $p\in [0,1]$, we also note
\begin{equation}
\label{bordaleatoire}
\partial_{p} A (\omega)=\{y\in\Zd \backslash A: \; \exists x\in A\quad \{x,y\}\in \mathcal G_p(\omega)\}.
\end{equation}

On this probability space, we now define a homogeneous Markov chain $(X_t)_{t \ge 0}$ with values in $S^{\Zd}$ and with transition probabilities as in~(\ref{matrix}).

\subsubsection*{Definition of the process}
Let $\xi^0 \in S^{\Zd}$ be a fixed initial configuration. 
We define 
$$
\begin{array}{lll}
A^{\blue}_0=\{x \in \Zd: \; \xi^0_x \in \{ {\blue},{\green}\}\}
& \text{and} & A^{\yellow}_0=\{x \in \Zd: \; \xi^0_x \in \{ {\yellow},{\green}\}\}, \\
B^{\blue}_0=\{x \in \Zd: \; \xi^0_x \in \{ \blue,\green,
{\blue}^*,{\green}^*\}\}
& \text{and} & B^{\yellow}_0=\{x \in \Zd: \; \xi^0_x \in \{ \yellow,\green,{\yellow}^*,{\green}^*\}\}.
\end{array}
$$
Note that by intersection and difference, we can exactly recover through these four sets the whole configuration $\xi^0$.

Let $0\le p_{\yellow} \le p_{\blue}\le 1$, and consider a Bernoulli configuration $\omega \in \Omega$, which will give the evolution rules of the process.
An infection can only travel from an active site of the corresponding color to an empty site, via an edge which is $p$-open in $\omega$ for the parameter $p$ associated to this infection, \ie either $p_{\yellow}$ or $p_{\blue}$. As before, an  active green site is to be imagined as a superposition of an active yellow site and an active blue site. So, if $e$ is an edge between an active green site and an empty site, then three cases arise:
if $0\le \omega_e \le p_{\yellow}$ then $e$ will transmit to $x$ both infections, 
if $p_{\yellow} \le \omega_e \le p_{\blue}$ then $e$ will only transmit to $x$ the blue infection,
while if $p_{\blue} \le \omega_e$ then no infection will travel through $e$ to $x$.

To determine the state at time $t+1$ of an empty site $x$ at time $t$, we look simultaneously at all edges between active sites at time $n$ and $x$: 
if all these edges transmit the same infection -- blue or yellow -- then $x$ takes this color and becomes active,
if these edges transmit infections of the two different types, then $x$ becomes green and active, 
and otherwise, $x$ remains empty.
Active sites at time $t$ become passive sites of the same color at time $t+1$.
These rules are translated in the following recursive definitions:
$$\left\{
\begin{array}{rcl}
A^{\yellow}_{t+1}(\omega) 
& = & \partial_{p_{\yellow}}A^{\yellow}_t(\omega) \backslash (B^{\blue}_t (\omega)\cup B^{\yellow}_t(\omega)), \\
B^{\yellow}_{t+1} (\omega)
& = & B^{\yellow}_t(\omega)\cup A^{\yellow}_{t+1}(\omega)
= B^{\yellow}_t (\omega)\cup (\partial_{p_{\yellow}}B^{\yellow}_t(\omega) \backslash
 B^{\blue}_t(\omega)), \\
A^{\blue}_{t+1} (\omega)
& = & \partial_{p_{\blue}}A^{\blue}_t(\omega) \backslash (B^{\blue}_t (\omega)\cup B^{\yellow}_t(\omega)), \\
B^{\blue}_{t+1} (\omega)
& = & B^{\blue}_t\cup A^{\blue}_{t+1}(\omega))
 =  B^{\blue}_t (\omega)\cup (\partial_{p_{\blue}}B^{\blue}_t(\omega) \backslash
 B^{\yellow}_t(\omega)).
\end{array}
\right.
$$
The set $A^{\blue}_t$ (\resp $A^{\yellow}_t$) is the set of active sites at time $t$ in that are either blue or green (\resp yellow or green), while $B^{\blue}_t$ (\resp $B^{\yellow}_t$) is the set of sites at time $t$ that are either blue or green (\resp yellow or green).
Note that by these definitions, a given site can be active at one time at most. 
We define then, for every $t\ge 0$, the value of the process $X_t$ at time $t$ as the element of $S^{\Zd}$ encoded by the four random sets $A^{\yellow}_t$, $B^{\yellow}_t$, $A^{\blue}_t$ and $B^{\blue}_t$.

\begin{lemme} \label{randomset}
The process $(X_{t})_{t\ge 0}$ is a homogeneous Markov chain governed by the transition probabilities defined in (\ref{matrix}).
\end{lemme}
\begin{proof}
Fix $\xi_0 \in S^{\Zd}$, and define $A^{\blue}_0$, $A^{\yellow}_0$, $B^{\blue}_0$ and $B^{\yellow}_0$ as previously.
 
The only point is to prove that $(X_{t})_{t\ge 0}$ is a homogeneous Markov chain, the identification of the transition probabilities is clear by construction.
The ideas of the proof stay in the following easy remarks:
\begin{itemize}
\item During the process, any site can only be active at one time at most.
\item Suppose that at time $t$, the process is in state $\xi$. To decide in which state it will switch at time $t+1$, the only edges that are to be examined are the ones between an active site and an empty site in $\xi$.
\item Thus, during the process, each edge is examined only once at most.
\end{itemize}
So knowing the present, the past will not affect the future.
In the rest of the proof, we try to turn this crude argument into a more rigorous one. 

In order to define the four random sets at time $t+1$ from the four random sets at time $t$ and  $\omega$, we introduce, 
for any subsets $A,B,C$ of $\Z^d$, any Bernoulli configuration $\omega \in \Omega$ and any probability $0\le p\le 1$, define the two following functions:
\begin{eqnarray*}
F(p,\omega,A,B,C) & = & \partial_{p}A(\omega) \backslash (B\cup C), \\
G(p,\omega,A,B) & = & A\cup (\partial_{p}A(\omega) \backslash B). 
\end{eqnarray*}
Then, the previous definitions are equivalent to:
$$
\left\{
\begin{array}{rcl}
A^{\yellow}_{t+1}(\omega) 
& = &  F(p_{\yellow},\omega,A^{\yellow}_t(\omega) ,B^{\yellow}_t(\omega) ,B^{\blue}_t(\omega) ), \\
B^{\yellow}_{t+1}(\omega)  
& = & G(p_{\yellow},\omega,B^{\yellow}_t(\omega) ,B^{\blue}_t(\omega) ), \\
A^{\blue}_{t+1} (\omega) 
& = & F(p_{\blue},\omega,A^{\blue}_t(\omega) ,B^{\yellow}_t(\omega) ,B^{\blue}_t(\omega) ), \\
B^{\blue}_{t+1} (\omega) 
& = & G(p_{\blue},\omega,B^{\blue}_t(\omega) ,B^{\yellow}_t(\omega) ),
\end{array}
\right.
$$
which is equivalent to say that $(X_t)_{t\ge 0}$ satisfies  
a recurrence formula of the type $X_{t+1}=f(X_t,\omega)$, where the function $f$ can be expressed in terms of the two functions $F$ and $G$.

To obtain the canonical Markov Chain representation $X_{t+1}=f(X_t,\omega^{t+1})$, we are going to build a coupling between
a random variable uniformly distributed on $\Omega$ and an independent and identically distributed sequence $(\omega^t)_{t\ge 1}$ with the same law.
Let $(\tilde{\Omega}, \tilde{\mathcal{F}},\tilde{\P})$ be a probability space and let $(\omega^t)_{t \ge 1}$ be independent $[0,1]^{\Ed}$-valued random variables with $\text{Unif}([0,1])^{\otimes\Ed}$ as common law.
We define $A^{\blue}_0$, $A^{\yellow}_0$, $B^{\blue}_0$ and $B^{\yellow}_0$ exactly as previously. But now, we set recursively:
\begin{eqnarray*}
A^{\yellow}_{t+1} = F(p_{\yellow},\omega^{t+1},A^{\yellow}_t,B^{\yellow}_t,B^{\blue}_t)
& \text{and} & 
B^{\yellow}_{t+1} = G(p_{\yellow},\omega^{t+1},B^{\yellow}_t,B^{\blue}_t), \\
A^{\blue}_{t+1} = F(p_{\blue},\omega^{t+1},A^{\blue}_t,B^{\yellow}_t,B^{\blue}_t)
& \text{and} & 
B^{\blue}_{t+1} = G(p_{\blue},\omega^{t+1},B^{\blue}_t,B^{\yellow}_t),
\end{eqnarray*}
Note that these four sets are measurable with respect to the $\sigma$-algebra generated by $(\omega^1,\dots,\omega^{t+1})$.
Let $\tilde{\omega}^0$ be a random variable defined on $(\tilde{\Omega}, \tilde{\mathcal{F}},\tilde{\P})$, with law $\text{Unif}([0,1])^{\otimes\Ed}$, and independent of the sequence $(\omega^t)_{t \ge 1}$, and define $(\tilde{\omega}^t)_{t \ge 0}$  recursively as follows: for any edge $e=\{x,y\}\in \Ed$, set  
\begin{eqnarray*}
\tilde{\omega}_e^{t+1} = &
\begin{cases}
\omega_e^{t+1} & \text{if }
x \in (A^{\yellow}_t \cup A^{\blue}_t) \text{ and }y \notin (B^{\yellow}_t \cup B^{\blue}_t),\\
\tilde{\omega}_e^t  & \text{otherwise.}
\end{cases}
\end{eqnarray*}
By natural induction, we prove that the law of $\tilde{\omega}^t$ under $\tilde{\P}$ is $\text{Unif}([0,1])^{\otimes\Ed}$.
By construction, each edge $e$ writes  $e=\{x,y\}$ with $x \in (A^1_t \cup A^2_t) \text{ and } y \notin (B^1_t \cup B^2_t)$ for at most one value of $t$.
It follows that the sequence $(\tilde{\omega}^t)_{t \ge 0}$ converges in the product topology to a limit that we denote $\tilde{\omega}^{\infty}$. Since the law of $\tilde{\omega}^t$ under $\tilde{\P}$ is $\text{Unif}([0,1])^{\otimes\Ed}$, it follows that the law of $\tilde{\omega}^{\infty}$ under $\tilde{\P}$ is also $\text{Unif}([0,1])^{\otimes\Ed}$.

Now, it is not difficult to see that the sequence $(X_t)_{t\ge 0}$ defined from $\tilde{\omega}^{\infty}$ as previously, satisfies the recurrence 
formula $X_{t+1}=f(X_t,\tilde{\omega}^{\infty})$, but also  $X_{t+1}=f(X_t,\omega^{t+1})$, which proves that $(X_{t})_{t\ge 0}$ is an homogeneous Markov chain.
\end{proof}

\subsection{Monotonicity properties and notations}

>From now on, we will denote by $(X_{t}^{\xi,p,q})_{t\ge 0}$ the competition process where
\begin{itemize}
\item $\xi \in S^{(\Zd)}$ is the initial configuration: $X_{0}^{\xi,p,q}=\xi$,
\item $0 \le p \le q \le 1$: the weakest (also called yellow) infection uses parameter $p$ while the strongest (also called blue) uses $q$.
\end{itemize} 
The corresponding random sets are now denoted by:
$$
\begin{cases}
\eta^{1}_{\xi,p,q}(t)=\{x\in\Zd: \; X^{\xi,p,q}_t(x)\in
\{\yellow,\yellowstar,\green,\greenstar\}\}, \\
\eta^{2}_{\xi,p,q}(t)=\{x\in\Zd: \; X^{\xi,p,q}_t(x)\in
\{\blue,\bluestar,\green,\greenstar\}\}.
\end{cases}
$$
Thus for $t\ge1$, they are also defined by the following recursive rules -- remember that the notation $\partial_p$ was defined in~(\ref{bordaleatoire}):
\begin{equation}
\label{THEdefinition}
\begin{cases}
\eta^{1}_{\xi,p,q}(t)= \eta^{1}_{\xi,p,q}(t-1)\cup (\partial_{p} \eta^{1}_{\xi,p,q}(t-1) \backslash \eta^{2}_{\xi,p,q}(t-1)) , \\
\eta^{2}_{\xi,p,q}(t)= \eta^{2}_{\xi,p,q}(t-1)\cup (\partial_{q} \eta^{1}_{\xi,p,q}(t-1) \backslash \eta^{2}_{\xi,p,q}(t-1)).
\end{cases}
\end{equation}
This particular realization of our competition process will be used in the sequel of the paper, because it presents the advantage to give an easy access to coupling and monotonicity properties.
Note that the function 
$$
\begin{array}{rrcl}
G: & [0,1] \times \Omega \times \mathcal{P}(\Zd) \times \mathcal{P}(\Zd) & \longrightarrow & \mathcal{P}(\Zd) \\
  & (p,\omega,A,B) & \longmapsto & G(p,\omega,A,B)=A\cup (\partial_{p}A(\omega) \backslash B)
\end{array}
$$
introduced in the proof of Lemme~\ref{randomset} is non-decreasing in $p$ and $A$, and non-increasing in $B$. As it defines the random sets at time $t+1$  from the random sets at time $t$, this implies in particular that:
\begin{lemme} $\;$ 
\label{yal}
\begin{itemize}
\item
$\eta^{1}_{\xi,p,q}(t+1)$ is non-decreasing in $p$ and non-increasing in $q$,
\item
$\eta^{2}_{\xi,p,q}(t+1)$ is non-decreasing in $q$ and non-increasing in $p$.
\end{itemize}
\end{lemme}

The next Lemma is trivial, but it is an illustration of the fundamental role played by the chemical distance in Bernoulli percolation in our analysis of this competition model: it says that the set of sites infected by any of the two infections at time $n$ can be compared with the single weaker infection.
\begin{lemme}
\label{comparaison}
Let us define, for any $0\le p\le 1$ and any $s \in \Zd$, the process $(B^s_p(t))_{t\in\N}$ by :
$$
B^s_p(0)=\{s\}
\text{ and }  \forall t \ge0, \;
B^s_p(t+1)=B^s_p(t) \cup \partial_{p} B^s_p(t).
$$
Let $s_1$ and $s_2$ be two distinct sites of $\Zd$ and 
$\xi$ be the element of $S^{\Zd}$ where all sites are empty, but $\xi_{s_1}=\yellow$ and $\xi_{s_2}=\blue$. Suppose that $0\le p\le q \le 1$. Then
$$\forall t\in\N\quad 
\left\{
\begin{array}{l}
B^{s_1}_p(t)\subset \eta^1_{\xi,p,q}(t)\cup \eta^2_{\xi,p,q}(t), \\
\eta^1_{\xi,p,q}(t) \subset B^{s_1}_p(t) \text{ and } 
\eta^2_{\xi,p,q}(t) \subset B^{s_2}_q(t).
\end{array}
\right.$$
\end{lemme}
It is easy to see that $B^s_p(t)=\{x\in\Zd: \;  D_p(s,x)\le t\},$
where $D_p(x,y)$ is the cardinal of the shortest $p$-open path from $x$ to $y$ and is called the \emph{chemical distance} between $x$ and $y$. Note that the inclusion $\eta^1_{\xi,p,q}(t) \subset B^{s_1}_p(t)$ implies that if $p<p_c$, then the infection with parameter $p$ almost surely dies out. 

This description~(\ref{THEdefinition}) of the competition model leads us to recall notations and results about chemical distance in Bernoulli percolation.

\section{Chemical distance in Bernoulli percolation}
\label{Sectionbernoulli}

In this section, we recall results concerning chemical distance in supercritical Bernoulli percolation:
\begin{itemize}
\item almost-sure convergence results~(\ref{asymptotic-speed}) and~(\ref{asymptotic-shape}) of the chemical distance to a deterministic norm,
\item large deviations inequalities~(\ref{ap-enrhume}) and~(\ref{shapeGD}) associated to this convergence,
\item classical estimates~(\ref{amasfini}) and~(\ref{amasinfini}) on the geometry of clusters.
\end{itemize}

We first complete the notations introduced at the beginning of Subsection~\ref{partbernouill}: the connected component of the site $x$ in the random graph $\mathcal{G}_p$ is denoted $C_p^x$, and the event that two sites $x$ and $y$ are in the same connected component of this graph is denoted $x \communiquep y$. 
Bernoulli percolation is in particular famous for its phase transition: there exists $0<p_c=p_c(d)<1$ such that
\begin{itemize}
\item if $p < p_c$ then with probability $1$, the random graph $\mathcal G_p$ has only finite connected components,
\item if $p > p_c$ then with probability $1$, the random graph $\mathcal G_p$ has at least one infinite connected component, which is moreover almost surely unique and denoted $C_p^\infty$.
\end{itemize}
See the reference book by Grimmett~\cite{Grimmett-book} for instance.

A \textit{path} is a sequence $\gamma=(x_1,
e_1,x_2,e_2,\ldots,x_n,e_n,x_{n+1})$ such that $x_i$ and $x_{i+1}$ are
neighbors and $e_i$ is the edge between $x_i$ and $x_{i+1}$. 
We will also sometimes describe $\gamma$ only by the vertices it
visits $\gamma=(x_1,x_2,\ldots,x_n,x_{n+1})$ or by its edges
$\gamma=(e_1,e_2,\ldots,e_n)$. The number $n$ of edges in $\gamma$ is called
the \textit{length} of $\gamma$ and is denoted by $|\gamma|$. A path is said to be $p$-\textit{open} in the configuration $\omega$
if all its edges are $p$-open in $\omega$.
The \emph{chemical distance} $D_p$ is the usual graph distance in $\mathcal G_p$:
$$\forall x,y \in \Zd \quad D_p(x,y)=\inf\{|\gamma|: \; \gamma \text{ $p$-open path between $x$ and $y$}\}.$$
We also define the random balls associated to this random distance:
$$\forall x \in \Zd, \, \forall t \ge 0 \quad 
B_p^x(t)=\{y \in \Zd: \quad D_p(x,y)\le t \}.$$
The formulation in terms of random distance comes from classical first-passage percolation, and indeed, this model can be seen as \iid first-passage percolation, where the passage-time of an edge takes value $1$ with probability $p$ and value $\infty$ with probability $1-p$.
An asymptotic shape result is also available for this model: in a previous paper~\cite{GM-fpppc}, we proved the existence of a deterministic norm $\|.\|_p$ on $\Rd$ such that $B^0_p(t)/t$ converges to the unit ball for $\|.\|_p$ on the event $\{0\communiquep\infty\}=\{0 \in C_p^\infty\}$, for the Hausdorff distance between two non empty compact subsets of $\Rd$.
For $x\in\Rd$ and $t\ge 0$, first define the deterministic balls associated to the norm $\|.\|_p$:
$$\mathcal{B}_{p}^x(t)=\{y\in\Rd: \; \|x-y\|_p\le t\}.$$
The Hausdorff distance between two non empty compact subsets  $A$ and 
$B$ of $\Rd$
is defined by
$$\mathcal{D}(A,B)=\inf\{t\ge 0: \; A\subset B+\mathcal{B}_{p}^0(t)\text{ and } 
B\subset A+\mathcal{B}_{p}^0(t)\}.$$
Note that the equivalence 
of norms on $\Rd$ ensures that the topology induced by this Hausdorff distance 
does not depend on the choice of the norm $\|.\|_p$. 
The convergence result writes then: for every $p>p_c(d)$,
\begin{itemize}
\item Existence of an asymptotic speed (Lemma 5.7 in~\cite{GM-fpppc}). 
\begin{equation}
\label{asymptotic-speed}
\miniop{}{\limsup}{\|x\|_p \to \infty} \1_{\{0 \communiquep x\}} \left( \frac{D_p(0,x)}{\|x\|_p} -1\right) =0 \quad \P \as
\end{equation}
\item Asymptotic shape result (Theorem 5.3 and Corollary 5.4 in~\cite{GM-fpppc}). \\
If $\Pcond_p(A)=\P(A| 0\communiquep\infty)$, then
\begin{equation}
\label{asymptotic-shape}
 \lim_{t\to +\infty}\mathcal{D}\left(\frac{B_p^0(t)}t,\mathcal{B}_{p}^0(1)\right)=0\quad\Pcond_p\as
\end{equation}
\end{itemize}
In the sequel, we will also use a corollary of these results.
For $A\subset\Zd$, we denote
$$|A|_p=\sup\{\|x\|_p: \; x\in A\} \quad
\text{ and }
\quad |A|_{*,p}=\inf\{\|x\|_p: \; x\in C_{p}^\infty \backslash A\}.$$
\begin{lemme}
\label{shape}
Let $p>p_c(d)$. On the event $\{a\in C^{\infty}_p\}$, we have $\P$ almost surely:
$$\frac{|B^a_{p}(t)|_{p}}t\to 1 \quad 
\text{ and }
\quad \frac{|B^a_{p}(t)|_{*,p}}t\to 1.$$
\end{lemme}
\begin{proof}
The identities 
$\limsup_{t\to +\infty}\frac{|B^a_{p}(t)|_{p}}t= 1$
and
$\limsup_{t\to +\infty}\frac{|B^a_{p}(t)|_{*,p}}t\le 1$
obviously follows from~(\ref{asymptotic-shape}).
It remains to show that for each $\delta>0$,
$$\P \left( \frac{|B^a_{p}(t)|_{*,p}}t\le 1-\delta\ \io \right)=0.$$
Suppose $\frac{|B^a_{p}(t)|_{*,p}}t\le 1-\delta\ \io$:
there exists  sequences $(x_n)_{n\ge 1}$ and  $(t_n)_{n\ge 1}$,
with $x_n\in C^\infty_{p}$, $\|x_n\|_p\le t_n(1-\delta)$, $D_p(0,x_n)\ge t_n$
and $t_n\to +\infty$. The sequence 
$(x_n)_{n\ge 1}$ is necessary unbounded, otherwise there would exist
a limiting value $x$, with $D_p(0,x)=+\infty$ and $x \in C^\infty_p$, which
is not possible.
It follows that there exist infinitely many $x\in C^\infty_p$ with
$D_p(0,x)\ge (1+\delta)\|x\|_p$.
By~(\ref{asymptotic-speed}), this happens with a null probability.
\end{proof}

As a direct consequence of these convergence results and of the coupling identity
$$D_q(0,nx) \1_{\{0\communiquep nx\}} 
\le D_p(0,nx) \1_{\{0\communiquep nx\}},
$$
we obtain the natural large comparison between norms for different parameters. It will be improved in Section~\ref{Sectionstrictecomp} to prove Proposition~\ref{Introstrictecomp}. 

\begin{lemme}
\label{Lemmelargecomp}
If $p_c(d)<p \le q \le 1$, then for every $x \in \Rd$, $\|x\|_q \le \|x\|_p$.
\end{lemme}

In another paper~\cite{GM-large}, we gave further information on the speed of convergence by
establishing the following large deviation inequalities corresponding to the previous convergence results:
for every $p>p_c(d)$, for every $\epsilon>0$, we have:
\begin{itemize}
\item Directional large deviation result.
\begin{equation}
\label{ap-enrhume}
\miniop{}{\limsup}{\Vert x\Vert_1\to +\infty} 
\frac{1}{\|x\|_1} \ln \P \left( 0\communiquep
  x, \; \frac{D_p(0,x)}{\|x\|_p}\notin (1-\epsilon, 1+\epsilon)
\right)<0.
\end{equation} 
\item Shape large deviation result.
There exist two strictly positive constants $A$ and $B$ such that
\begin{equation}
\label{shapeGD}
\forall t>0 \quad \Pcond_p \left(
\mathcal{D}\left(\frac{B_p^0(t)}t,\mathcal{B}^0_{p}(1)\right) \ge \epsilon
\right) \le A e^{-Bt}.
\end{equation}
\end{itemize}

As a consequence, we obtain the next lemma, which enables the control of minimal paths:
\begin{lemme}
\label{chemin-pas-loin}
Note $H_p(x,y,\epsilon)$ the event: ``There exists a $p$-open minimal  path from
$x$ to $y$  
which is completely inside $\mathcal{B}_{p}^y((1+\epsilon)\|x-y\|_p)$
and whose length 
is smaller than $(1+\epsilon)\|x-y\|_p$''.

Then for every $p>p_c(d)$,
for every $\epsilon>0$, there exist two strictly positive constants $A$
and $B$ such that:
$$\forall x,y\in\Zd \quad \P(y\in C_p^\infty,\, x\in C_p^\infty, \;
H_p(x,y,\epsilon)^c)\le A \exp(-B \|x-y\|_1).$$
\end{lemme}
\begin{proof}
Using translation invariance, we can assume that $y=0$.
Note that $H(x,0,\epsilon)^c$ contains the event 
$\{D_p(0,x) \le \|x\|_p (1+\epsilon/2) \} \cap
\{B_p^0(t (1+ \epsilon/2)) \subset
\mathcal{B}_p^0((1+\epsilon)\|x\|_p)\}$ and apply 
the large deviation inequality for the chemical distance~(\ref{ap-enrhume}) and the large deviation inequality~(\ref{shapeGD}) for the asymptotic shape.
\end{proof}  

We also recall here some classical results concerning the geometry of clusters in supercritical percolation.
Thanks to Chayes, Chayes, Grimmett, Kesten and Schonmann~\cite{CCGKS}, we 
can control the radius of finite clusters: there exist two strictly positive 
constants $A$ and $A$ such that
\begin{equation}
\label{amasfini}
\forall r>0 \quad
 \P \left( |C^0_p|<+\infty, \; 0\communiquep\partial_1 \mathcal{B}_1^0(r) \right)
\le A e^{-B r}.
\end{equation}
The size of holes in the infinite cluster can also be controlled: there exist 
two strictly positive constants $A$ and $B$ such that
\begin{equation}
\label{amasinfini}
\forall r>0 \quad
 \P \left( C^\infty_p \cap \mathcal{B}_1^0(r) =\varnothing \right)
\le Ae^{-B r}.
\end{equation}
When $d=2$, this result follows from large deviation estimates by Durrett 
and Schonmann~\cite{DS}. Their methods can easily be transposed when $d\ge 3$. 
Nevertheless, when $d\ge 3$, the easiest way to obtain it seems to use Grimmett 
and Marstrand~\cite{Grimmett-Marstrand} slab's result.   

Note that in Lemma~\ref{chemin-pas-loin}, in $(\ref{amasfini})$ and in $(\ref{amasinfini})$, thanks to the norm equivalence, the choice of
the norm $\|.\|_1$ is of course irrelevant, but in the very values of the positive constants.

\section{Strict inclusion of asymptotic shapes for chemical distance}
\label{Sectionstrictecomp}

Inequalities on asymptotic shapes are already known for classical first-passage
percolation -- see the papers by Van den Berg and Kesten~\cite{vdb-kes} and by  Marchand~\cite{marchand}. The aim of this section is to prove Proposition~\ref{Introstrictecomp}, which is the analogous result in this context. We recall that the large inequality was easily established in Lemma~\ref{Lemmelargecomp}, but that strict comparisons will be crucial to handle the competition problem. 

The proof of Proposition~\ref{Introstrictecomp} is based on renormalization techniques. We thus begin by stating an adapted renormalization lemma, which is the one used by Van den Berg and Kesten in~\cite{vdb-kes}.

\subsection{A renormalization lemma.}

\subsubsection*{The renormalization grid}
\begin{figure}[h!]
\begin{flushleft}
\input{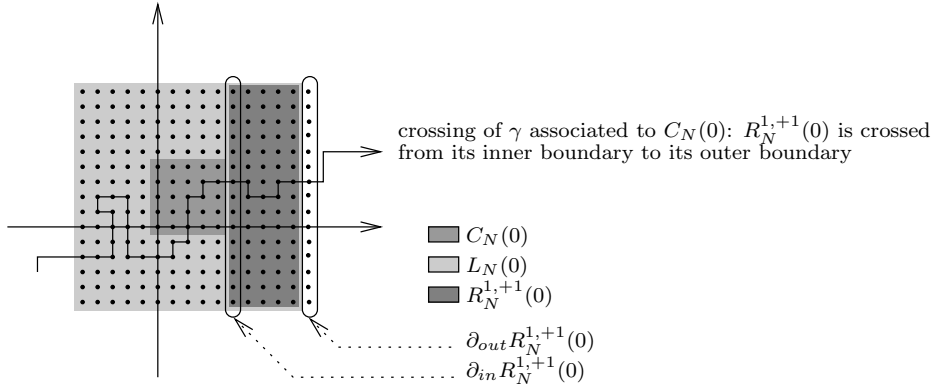}
\caption{Elements of the renormalization grid for $N=6$ in dimension $d=2$.} 
\label{DESSgrille}
\end{flushleft}
\end{figure}

Let $N$ be a strictly positive integer. We introduce the following
notations. 
\begin{itemize} 
\item
$C_N$ is the cube $[-1/2,N-1/2]^d$. We call $N$-{cubes} the
cubes $C_N(k)=kN+C_N$ obtained by translating $C_N$ according to $Nk$ with $k
\in \Zd$. The coordinates of $k$ are called the coordinates of the $N$-cube
$C_N(k)$. Note that $N$-cubes induce a partition of $\Zd$.  
\item
$L_N$ is the large cube $[-N-1/2,2N-1/2]^d$, and the large cube $L_N(k)$ is
obtained by translating $L_N$ according to $Nk$ with $k \in \Zd$. The
boundary of $L_N(k)$, denoted by $\partial L_N(k)$, is the set of sites outside
$L_N(k)$ that have a neighbor in $L_N(k)$. 
\item
$R_{N}$ is the rectangular box $[-1/2,N-1/2]^{d-1}
\times[-N-1/2,2N-1/2]$. In the large cube $L_N(k)$, the $N$-cube
$C_N(k)$ is surrounded by the $2d$ $N$-{boxes}, obtained by rotations and
translations of $R_N$. For instance, in $L_N(0)$, the $N$-cube $C_N(0)$ is surrounded by the $2d$ following $N$-boxes: for $1\le i \le d$, and for $\signe \in \{-1,+1\}$, we define
$$R_{N}^{i,\signe}(0)=\left[-N-\frac12,2N-\frac12 \right]^{i-1} 
\times \left[\signe N- \frac12 ,(1+\signe)N-\frac12 \right] \times 
\left[-N-\frac12,2N-\frac12\right]^{d-i}.$$
The set of all these surrounding boxes is denoted $\mathcal R_N$.
\end{itemize}
An edge is said to be \emph{in a subset $E$} of $\Rd$ if at least one of
its two extremities is in $E$.
We now define the \emph{inner} and \emph{outer} boundaries of a $N$-box associated
to a pair $(C_N(k),L_N(k))$ of cubes. Let's do this for $R_{N}^{1,+}(0)$
and extend the definition to other boxes by rotation and translation:
$$
\begin{array}{rcl}
\partial_{out}R_{N}^{1,+1}(0) & = & \{ (2N, y), \; y \in [-N, \dots,
2N-1]^{d-1}\}, \\
\partial_{in}R_{N}^{1,-1}(0) & = & \{ (N, y), \; y \in [-N, \dots, 2N-1]^{d-1}\}.
\end{array}
$$
Note that $\partial L_N(k)$ is the disjoint union of the sets
$(\partial_{out}R_{N}^{i,\signe}(k))_{1 \leq i \leq d, \; \signe\in \{+1,-1\}}$, and that a path entering in $C_N(k)$
and getting out of $L_N(k)$ has to cross one of the $2d$ $N$-boxes
surrounding $C_N(k)$ in $L_N(k)$, from its inner boundary to its outer
boundary. We can then define the crossing
associated to a $N$-cube $C_N(k)$ -- see also Figure~\ref{DESSgrille}:
\begin{defi} \label{DEFcrossing}
Let $\gamma=(x_0, \dots , x_l)$ be a path such that $x_0 \in C_N(k)$ and
$x_l \not\in L_N(k)$. We set $j_f=\min \{0\leq k \leq l, \; x_k \in
\partial L_N(k)\}$. There exists a unique $(i,\signe)$ such that $x_{j_f} \in
R_{N}^{i, \signe}(k)$. Let then $j_0=\max \{0 \leq k \leq j_f, \; x_k \notin
R_{N}^{i, \signe}(k)\}$. The portion $(x_{j_0+1}, \dots , x_{j_f})$ of $\gamma$ is the
crossing of $\gamma$ associated to $C_N(k)$.
\end{defi}

\subsubsection*{Main crossings of a path.}
Let $N$ be a strictly positive integer, $x$ be a point in $\Zd$ and
$\gamma$ be a path without any double point from $0$ to $x$.
We want to associate to $\gamma$ a sequence of crossings of $N$-boxes (the main
crossings of $\gamma$),
in a way that two different crossings are edge-disjoint. Consider
first the sequence $\sigma_0=(k_1,
\dots, k_{\tau_0})$ made of the coordinates of the $N$-cubes
successively visited by $\gamma$. As the $N$-cubes induce a partition of $\Zd$, this sequence is well defined, and has the following properties:
$$
(\mbox{\textnormal{P}}_0) \left\{
\begin{array}{l}
0 \in C_N(k_1), \; \;  x \in C_N(k_{\tau_0}), \\
\forall \; 1 \leq i \leq \tau_0-1, \; \|k_{i+1}-k_i\|_1=1.
\end{array}
\right.
$$
But $\sigma_0$ can have doubles points; we remove them by the
classical loop-removal process described in~\cite{grimmett-kesten}. We thus obtain a sequence $\sigma_1=(k_{\phi_1(1)},
\dots, k_{\phi_1(\tau_1)})$ extracted from $\sigma_0$ , 
with the following properties: 
$$
(\mbox{\textnormal{P}}_1) \left\{
\begin{array}{l}
0 \in C_N(k_{\phi_1(1)}), \; \; x \in C_N(k_{\phi_1(\tau_1)}), \\
\forall \; 1 \leq i \leq \tau_1-1, \; \|k_{\phi_1(i+1)}-k_{\phi_1(i)}\|_1=1,\\
\sigma_1 \mbox{ has no double point.}
\end{array}
\right.
$$
To every cube $C_N(k)$ in this sequence such that $\gamma$ gets out of
$L_N(k)$, that means for every $N$-cube in $\sigma_1$ with the possible exception of the $2d$ last, we associate a crossing of a
$N$-box in the following way:
let $z$ be the first point of $\gamma$ to be in
$C_N(k)$, and let $z_2$ be the first point of $\gamma$ after $z$ 
to be in $\partial L_N(k)$. Then the crossing associated to the
$N$-cube $C_N(k)$ is the crossing of the portion of $\gamma$ between
$z$ and $z_2$ associated to $C_N(k)$ in Definition~\ref{DEFcrossing}. 

The problem now is that two distinct cubes in $\sigma_1$
can have the same associated crossing. We have to extract a subsequence 
once again in order to obtain edge-disjoint crossings. Set $\phi_2(1)=1$, 
and define $\phi_2$ by induction:
$$\phi_2(i+1)= \inf \{ j>\phi_2(i) \mbox{ such that }
\|k_{\phi_1(j)}-k_{\phi_1 \circ \phi_2(i)}\|_\infty>1\} -1
$$
if the infimum exists, and let $\tau$ be the smallest index $i$ for
which $\phi_2(i+1)$ is not defined. Set
$\phi=\phi_1 \circ \phi_2$; the elements of 
$\sigma=(k_{\phi(i)})_{1 \leq i \leq \tau}$ are called the
{main cubes}
of $\gamma$, and their associated crossings the \emph{main crossings of}
$\gamma$. This sequence has the following properties (see~\cite{vdb-kes}):
$$
(P) \left\{
\begin{array}{l}
0 \in C_N(k_{\phi(1)}), \\
\|k_{\phi(\tau)}-k_{\tau_0}\|_\infty \leq 1,\\
\forall \; 1 \leq i \leq \tau-1, \; \|k_{\phi(i+1)}-k_{\phi(i)}\|_\infty =1,\\
\mbox{the main crossings of }\gamma \mbox{ are edge-disjoint.}
\end{array}
\right.
$$
>From Properties $(P)$ we can deduce that for every $x$ in $\Zd$, the number $\tau$ of main $N$-cubes 
of a path with no double point from $0$ to $x$ satisfies the following inequality:
\begin{equation} \label{LEMmintau}
\tau \geq \frac{\|x\|_\infty}{N}.
\end{equation}

\subsubsection*{A renormalization lemma.}
The following lemma is an adaptation of Lemma (5.2) in~\cite{vdb-kes}, 
and its proof is a standard Peierl's argument (see proof of (3.12) in~\cite{grimmett-kesten}).
We thus just state it without any proof.
\begin{lemme} \label{LEMren}
For each $N \in \N$, we give to the $N$-cubes
a random color, black or white, according to the states of the edges
in the initial model, such that:
\begin{itemize}
\item
For each $N \in \N^*$, the colors of the $N$-cubes are identically
distributed.
\item
For each $N \in \N^*$, for each $k \in \Zd$, the color of the $N$-cube
$C_N(k)$ depends only on the states of the edges in $L_N(k)$.
\item
$ \displaystyle \lim_{N \rightarrow + \infty} \P(C_N(k) \mbox{is black})=1$.
\end{itemize} 

\noindent
Then for every $\rho \in ]0,1[$, there exists $N_\rho$ such that
for all $N \geq N_\rho$, there exist two strictly positive 
constants $A$ and  $B$ such that for every $x \in \Rd$:
\begin{equation} \label{LEMrenEQU}
\P\left(
\begin{array}{c}
\mbox{There exists a path } \gamma \mbox{ from } 0 \mbox{ to } x
\mbox{ that, among }\\
\mbox{its } \tau \mbox{ main $N$-cubes, has less than } \rho \tau \mbox{ black cubes}
\end{array}
\right)   
\leq A \exp (-B \|x\|_\infty). 
\end{equation}
\end{lemme}

\subsection{Proof of the strict comparison result Proposition~\ref{Introstrictecomp}.} $\;$

Fix $p$ and $q$ such that $p_c(d)<p<\overrightarrow{p_c}(d)$ and $p<q \le 1$. Roughly speaking, as $p<\pcfleche(d)$, we can find along a $p$-minimal path from $0$ to $nx$ a certain number of crossing of rectangular boxes such that:
\begin{itemize}
\item the restriction of the $p$-minimal path of one box is not direct,
\item by adding $q$-open edges, as $q>p$, we can find in this box a direct
$q$-minimal path with the same extremities, which is thus an improvement of the $p$-minimal path.
\end{itemize}
By using these improvements, we can exhibit a significant discrepancy, \ie of order $n$, between
$D_p(0,nx)$ and $D_q(0,nx)$.
The proof consists in giving estimates to these crude arguments.

\begin{proof}
Consider the space
$\Omega=\{0,1\}^{\Ed}\times \{0,1\}^{\Ed},$ endowed with the classical Borel $\sigma$-algebra on $\Omega$ and the probability measure
$$\P=\Ber(p)^{\otimes\Ed}\otimes \Ber \left(\frac{q-p}{1-p} \right)^{\otimes\Ed}.$$
Write points of $\Omega$ in the following manner 
$$\omega=(\omega^1, \omega^2) \text{ with } \omega^1 =(\omega^1_e)_{e \in \Ed} \in \{0,1\}^{\Ed} \text{ and }\omega^2 =(\omega^2_e)_{e \in \Ed} \in \{0,1\}^{\Ed}.$$ Define then, for every $e \in  \Ed$, $\omega^3_e=\omega^1_e \vee \omega^2_e$. Clearly, the law of $(\omega^1_e)_{e\in\Ed}$ under $\P$ is $\Ber(p)^{\otimes\Ed}$ whereas the law of $(\omega^3_e)_{e\in\Ed}$ under $\P$ is $\Ber(q)^{\otimes\Ed}$. We denote by $\mathcal{G}_p$ -- \resp $\mathcal{G}_q$ --  the corresponding
random graphs and by $D_p(x,y)$ -- \resp $D_q(x,y)$ -- the random distance
from $x$ to $y$ in $\mathcal{G}_p$ -- \resp $\mathcal{G}_q$. Note that in this special coupling, 
$$\mathcal{G}_p \subset \mathcal{G}_q \text{ and } D_q(x,y) \le D_p(x,y).$$

For each $N \in \N$, we consider the same renormalization grid as previously
and give to each $N$-box ${R}_N^{i, \signe}(k)$ a random color:
\begin{defi}
The box ${R}_N^{i, \signe}(k)$, with $k \in \Zd$, is said to be black if and
only if it satisfies the following property:
\begin{equation*}
\begin{array}{c}
\forall y \in \partial_{in}R_N^1(k), \; 
\forall z \in \partial_{out}R_N^1(k), \; 
\forall \gamma \mbox{ $p$-open path from } y \mbox{ to } z 
\mbox{ included in  } R_N^1(k), \\ 
|\gamma|\geq \|z-y\|_1+1.
\end{array}
\end{equation*}
It is said to be white otherwise. This definition is naturally extended to other
boxes by translation and rotation.
\end{defi}
Thus a box is black if and only if it can not be directly crossed from its inner boundary to its outer boundary by a $p$-open path.
Let us verify that this coloring satisfies the conditions of
renormalization Lemma~\ref{LEMren}. It is clear that the colors of the different cubes are
identically distributed, and that the color of ${C}_N(k)$ only depends
on the states of the edges in $L_N(k)$. Let us now estimate
the probability $p_N$ for ${C}_N(0)$ to be white. 
It is clear by translation invariance that $p_N \leq 2d \P \left({R}_N^{1, +1}(0) \mbox{ is white}
 \right)$ and that the probability for ${R}_N^{1, +1}(0)$ to be white is
bounded by
$$
\P \left({R}_N^{1, +1}(0) \mbox{ is white}
 \right) \le (3N+1)^{d-1} 2^{d-1} \P \left( \max \{ \|x\|_1: \; x \in \overrightarrow{C_p}^0 \} \ge N \right),
$$
where $\overrightarrow{C_p}^0$ is the cluster containing $0$ in \emph{oriented} percolation with parameter $p$. The term $(3N+1)^{d-1}$ counts the possible starting points of the oriented open path, while the term $2^{d-1}$ counts its possible orientations. As in the non-oriented case, when $p <\overrightarrow{p_c}(d)$, the probability in the left-hand side member decreases exponentially fast with $N$ -- see the paper by Aizenman and Barsky~\cite{aizenman-barsky} --
which proves that:
\begin{equation*} 
\lim_{N \rightarrow
+\infty} p_N =0.
\end{equation*}
We can then apply the renormalization
Lemma~\ref{LEMren} with a fixed parameter $\rho$ satisfying
$0<\rho<1$. 
Let $N$
be large enough to have (\ref{LEMrenEQU}) with positive $A$ and $B$. These $\rho$ and $N$ are now fixed for the sequel of the proof.

For each $n \ge 1$ and every $x \in \Zd \backslash \{0\}$, if the event $\{0 \communiquep nx\}$ occurs, we denote by $\gamma_{n,x}$ be a $p$-open path from $0$ to $nx$ whose length is equal to $D_p(x,y)$.
Let $\sigma_{n,x}=(k_1, \dots, k_{\tau_{n,x}})$ be the sequence of its main
cubes and denote by $A_{n,x}$ the event that among
these $\tau_{n,x}$ main cubes, at most $\rho
\tau_{n,x}$ cubes are black. With Lemma~\ref{LEMren} we have:
\begin{eqnarray} 
\P(A_{n,x}\cap \{0\communiquep nx\}) & = & \P \left(
\begin{array}{c}
\mbox{there exists a $p$-open path } \gamma \mbox{ from } 0 \mbox{ to } nx \\
\mbox{that, among its } \tau \mbox{ main $N$-cubes,} \\
\mbox{has less than } \rho \tau \mbox{ black cubes}
\end{array} 
\right) \nonumber \\
& \leq  & A \exp (-B n\|x\|_\infty).
\label{LEMmajPAn}
\end{eqnarray}

We define now the notion of \emph{good rectangular boxes} 
\begin{defi}
A rectangular box $R$ is \emph{good} if it is black and  if, moreover, for every $e \in R$, $\omega^3(e)=1$.
\end{defi}
In other words, in a good box, edges that are not $p$-open are $q$-open.
Let $n$ be large enough and let $R\in \mathcal R_N$ be a good box. Suppose that the path $\gamma_{n,x}$ crosses $R$
 and that this crossing,
denoted by ${\gamma_{n,x}}_{|R}$, is a main crossing of $\gamma_{n,x}$.  Denote
by $y$ and $z$ the extremities of the restriction ${\gamma_{n,x}}_{|R}$ of
the path $\gamma_{n,x}$ to the box $R$. Then, by definition of \emph{black} and \emph{good} boxes, 
\begin{equation}
\label{cestmieux}
\|z-y\|_1=D_q(y,z) \le D_p(y,z)-1.
\end{equation}
Note that moreover, in this case, any $q$-open path between $y$ and $z$ with length $\|z-y\|_1=D_q(y,z)$ is completely inside $R$. Choose one and call it an \emph{improvement} for $D_q$ of  $\gamma_{n,x}$ in $R$.

Now, on the event $\{0 \communiquep nx\}$, replace in $\gamma_{n,x}$ all the restrictions 
associated to main crossings of $\gamma_{n,x}$ by their improvements for
$D_q$, to obtain a modified path $\hat \gamma_{n,x}$ from $0$ to $nx$:
this is possible, because by definition, main crossings are in
non-intersecting boxes. Then
\begin{eqnarray*}
&& \1_{\{0\communiquep nx\}} (D_p(0,nx)-D_q(0,nx)) \\
& \ge & \1_{\{0\communiquep nx\}}
\sum_{R \in \mathcal R_N} 
\1_{\{R \mbox{\scriptsize{ is good}}\}}
\1_{
\left\{
\begin{subarray}{c}
\gamma_{n,x} \text{ crosses } R, \text{ and} \\
\text{this is a main crossing of } \gamma_{n,x}
\end{subarray}
\right\}
}\\
& \ge & \1_{\{0\communiquep nx\}}
\sum_{R \in \mathcal R_N} 
\left(
\prod_{e \in R}
\1_{\omega^3_e=1}
\right) 
\1_{\{R \mbox{\scriptsize{ is black}}\}}
\1_{
\left\{
\begin{subarray}{c}
\gamma_{n,x} \text{ crosses } R, \text{ and} \\
\text{this is a main crossing of } \gamma_{n,x}
\end{subarray}
\right\}
} \\
& \ge & \1_{\{0\communiquep nx\}}
\sum_{R \in \mathcal R_N} 
\left(
\prod_{e \in R}
\1_{\omega^2_e=1}
\right) 
\1_{\{R \mbox{\scriptsize{ is black}}\}}
\1_{
\left\{
\begin{subarray}{c}
\gamma_{n,x} \text{ crosses } R, \text{ and} \\
\text{this is a main crossing of } \gamma_{n,x}
\end{subarray}
\right\}
}.
\end{eqnarray*}
Note $G(R)$ the event $\{\forall e \in R \quad \omega^2_e=1\}$.
As $\omega_1$ and $\omega_2$ are independent, the conditional law of the random variable
$$Y_{n,x} \stackrel{\text{def}}{=} \sum_{R \in \mathcal R_N} 
\1_{G(R)} \1_{\{R \mbox{\scriptsize{ is black}}\}}
\1_{
\left\{
\begin{subarray}{c}
\gamma_{n,x} \text{ crosses } R, \text{ and} \\
\text{this is a main crossing of } \gamma_{n,x}
\end{subarray}
\right\}
}
$$
knowing $\omega^1$ is a binomial law $\text{Bin}(Z_{n,x},r)$ with parameters
\begin{eqnarray*}
Z_{n,x} & \stackrel{\text{def}}{=} & \sum_{R \in \mathcal R_N} 
\1_{\{R \mbox{\scriptsize{ is black}}\}}
\1_{
\left\{
\begin{subarray}{c}
\gamma_{n,x} \text{ crosses } R, \text{ and} \\
\text{this is a main crossing of } \gamma_{n,x}
\end{subarray}
\right\}
}, \\
r & = & \P \left( G(R)\right) \ge \left( \frac{q-p}{1-p} \right) ^{c_d N^d} >0.
\end{eqnarray*}
We have then, using Estimate~(\ref{LEMmintau}) on the event $A^c_{n,x}$:
\begin{equation*}
{\1}_{A^c_{n,x}} 
{\1}_{\{0\communiquep nx\}} Z_{n,x}
\geq  \rho \tau_{n,x} {\1}_{A^c_{n,x}}{\1}_{\{0\communiquep nx\}}
\geq   \rho \frac{n\|x\|_\infty}{N} {\1}_{A^c_{n,x}}{\1}_{\{0\communiquep nx\}}. 
\end{equation*}
Thus, if $\delta>0$, we have
\begin{eqnarray*}
&& \P \left( 0\communiquep nx,\; 
Y_{n,x} \le \frac{ \rho n\|x\|_{\infty}}N r(1-\delta) \right) \\
& \le & \P(\{0\communiquep nx\} \cap A_{n,x})+
\sum_{k=\frac{\rho n\|x\|_{\infty}}N}^\infty \P \left( Z_{n,x}=k, \; Y_{n,x} \le \frac{\rho n\|x\|_{\infty}}N r (1-\delta) \right), \\
& \le & \P(\{0\communiquep nx\} \cap A_{n,x})+
\sum_{k \ge \frac{\rho n\|x\|_{\infty}}N} \P \left( Z_{n,x}=k, \; Y_{n,x} \le k r (1-\delta) \right), \\
& \le & \P(\{0\communiquep nx\} \cap A_{n,x})+
\sum_{k \ge \frac{\rho n\|x\|_{\infty}}N} \P \left( Z_{n,x}=k \right)
2 \exp\left( -\frac{k \delta^2}{4r(1-r)}\right) \\
&& \; \; \; \; \; \; \; \; \; \text{by Chernov inequality},\\
& \le & \P(\{0\communiquep nx\} \cap A_{n,x})+
2 \exp\left( -\frac{ \rho n\|x\|_{\infty}\delta^2}{4Nr(1-r)}\right).
\end{eqnarray*}
By (\ref{LEMmajPAn}) and  
Borel-Cantelli Lemma, this leads to
\begin{eqnarray*}
&& \P \left( 0\communiquep nx, \;  \frac{D_p(0,nx)}n-\frac{D_q(0,nx)}n\le \frac{\rho\|x\|_{\infty}(1-\delta)r}N\ \io \right)=0.
\end{eqnarray*}
On the event $\{0\communiquep \infty\} \subset \{0\communiqueq \infty\}$, by the convergence result~(\ref{asymptotic-speed}), we obtain 
$\| x \|_p-\| x \|_q\ge \frac{\rho\| x \|_{\infty}(1-\delta)r}{N}$, 
and finally, by letting $\delta$ going to $0$,
$$\forall x \in \Zd \quad \| x \|_p-\| x\|_q\ge \frac{\rho r}{N}\| x\|_{\infty}.$$
Since norms are homogeneous and continuous, this ends the proof.
\end{proof}

\section{Coexistence can only happen at slow speed}
\label{Sectionslowspeed}

We tackle in this section the core of the paper: the proof of Proposition~\ref{Introslowspeed}. 
For $p$ and $q$ larger than $p_c$, we define
$$C_{p,q}= \sup_{x \in \Rd \backslash \{0\}} \frac{\|x\|_q}{\|x\|_p}.$$
We fix here $p_1<p_2$ and two distinct sites $s_1$ and $s_2$ of $\Zd$: the initial state $\xi$ is the configuration where every site is empty, but $s_1$, which  is active yellow, and $s_2$, which is active blue. In the sequel, to lighten notations, we omit the subscripts $p_1,p_2,\xi$: for instance, 
$$\eta^2(t)=\eta^2_{p_1,p_2,\xi}(t).$$
By Proposition~\ref{Introstrictecomp}, we know that 
$C_{p_1,p_2}<1$. 

In fact, Proposition~\ref{Introslowspeed} will appear as a by-product of the
following theorem, which  ensures that if the $p_1$-infection survives, 
then the time of infection of $x$ by the $p_2$-infection, when it is
finite, should be of order $\|x\|_{p_1}$ rather than $\|x\|_{p_2}$, expected time of infection for one simple $p_2$-infection.

Define, for $x\in\Zd$ :  
\begin{eqnarray*}
&& t(x)  =  \inf \{ t\ge 0: \;  x\in\eta^2(t) \}, \\
&& \mathcal{G}^i =  \left\{ \sup_{t \ge 0} |\eta^i(t)|=+\infty \right\} \text{ for } i=1,2.
\end{eqnarray*}

\begin{theorem}
\label{avantage2}
Let $\delta>0$.
Then there exist $A,B>0$ such that 
$$
\forall x \in \Zd \quad \P \left( \mathcal{G}^1\cap \{t(x) \le (1-\delta)
\|x\|_{p_1} \} \right)\le A\exp(-B\|x\|).
$$
\end{theorem}

At first, let us see how Theorem~\ref{avantage2} implies Proposition~\ref{Introslowspeed}: 

\begin{proof}
Let $\delta>0$. We must prove that
$\P(\mathcal{G}^1, \;|\eta^2(t)|_{p_1}\ge (1+\delta)t \io)=0$.
Obviously, it is equivalent to prove that
$$\P(\mathcal{G}^1,\;  (1+\delta)t(x)\le \|x\|_{p_1} \text{ for infinitely many }x)=0.$$
This comes from Theorem~\ref{avantage2}, with the help of Borel-Cantelli's lemma.
\end{proof}

We still need some extra notations and lemmas.\\

\noindent\textbf{Definitions.}\\
We note $\mathcal{S}=\{x\in\Rd\ : \|x\|_{p_2}=1\}$ and define
the \emph{shells}:
for each $A \subset \mathcal{S}$, and every $0<r<R$, we set
\begin{eqnarray*}
\hat x & = & x/ \|x\|_{p_2}, \\ 
\shell(A,r,R) & = & \{x \in \Zd: \; \hat x \in A \text{ and } r \le \|x\|_{p_2}
\le R\}.
\end{eqnarray*}
So roughly speaking, $A$ is to think about as the set of possible
directions for the points in the shell, while $[r,R]$ is the set of radii.

For $A \subset \mathcal{S}$ and $\phi>0$, define the following
enlargement of $A$:
\begin{eqnarray*}
A \oplus \phi & = & (A+\mathcal{B}^0_{p_2}(\phi))\cap \mathcal{S}.
\end{eqnarray*}

\begin{lemme}
\label{direction}
For any norm $|.|$ on $\Rd$, one has
$$
\forall x,y \in \Rd \backslash\{0\} \quad 
\left| \frac{x}{|x|}-\frac{y}{|y|} \right|  
\le\frac{2|x-y|}{\max\{|x|,|y|\}}.
$$  
\end{lemme}

\begin{lemme}
\label{contunif}
For every $\rho>0$, there
exists $\theta>0$ such that
$$ \forall x,y \in \Rd \backslash\{0\} \quad 
\|\hat x -\hat y\|_{p_2} \le \theta \; \Longrightarrow \;(1-\rho)
\frac{\|x\|_{p_1}}{\|x\|_{p_2}} \le \frac{\|y\|_{p_1}}{\|y\|_{p_2}} \le (1 +\rho) \frac{\|x\|_{p_1}}{\|x\|_{p_2}}.$$
\end{lemme}

\begin{proof}
Note $F(x)=\frac{\|x\|_{p_1}}{\|x\|_{p_2}}$. Then, $\left|\frac{F(y)}{F(x)}-1\right|=\frac{|F(y)-F(x)|}{F(x)}\le C_{p_1,p_2}|F(x)-F(y)|$.\\
Now we have
$$
|F(x)-F(y)| = | F(\hat{x})-F(\hat{y})|
 = | \,\| \hat{x}\|_{p_1} -\|\hat{y}\|_{p_1} |
\le \| \hat{x}-\hat{y}\|_{p_1}
 \le  C_{p_2,p_1}\| \hat{x}-\hat{y}\|_{p_2}.
$$
Thus, we can take $\theta=\frac{\rho}{C_{p_1,p_2}C_{p_2,p_1}}>0$.
\end{proof}

We can now begin the proof of Theorem~\ref{avantage2}, which is cut into three main steps.

\subsection{Initialization of the spread} 

The aim of the next lemma is to see that if the event $\{ t(x) \le
(1-\delta)\|x\|_{p_1} \}$ is realized, then with high probability, at the
slightly largest time $(1-\delta')\|x\|_{p_1}$, the $p_2$-infection has
colonized a small shell, and this will provide it a strategic advantage for the next steps of the spread. 

\begin{lemme} 
\label{petit-bout-de-coquille}  
Let $\delta>0$ and choose any $0<\delta'<\delta$. 

For any $x \in \Zd \backslash \{0\}$, any $1<\gamma<\gamma'$ and any
$\theta>0$, we define the following events, depending on $x,\gamma,
\gamma'$ and $\theta$:
\begin{eqnarray*}
E_0 & = &  \{x\in C_{p_2}^\infty\}, \\ 
E_1 & = & \{\eta^1( (1-\delta')\|x\|_{p_1}) \subset
\mathcal{B}_{p_1}^0(\|x\|_{p_1}) \}, \\
E_2 & = & \{\eta^2((1-\delta')\|x\|_{p_1}) \, \supset \, C_{p_2}^\infty \cap
 \shell(\{\hat x\}\oplus\theta, \gamma\|x\|_{p_2},\gamma'\|x\|_{p_2})  \}, \\ 
E & = & E_0 \cap E_1 \cap E_2.
\end{eqnarray*}

Then there exist $\gamma'_0>1$ and $\theta_0>0$ such that for any
$1<\gamma<\gamma'<\gamma'_0$ and any $0<\theta<\theta_0$, there exist two strictly positive constants $A$ and $B$ such that 
$$ 
\forall x \in \Zd \quad \P( \{ t(x) \le (1-\delta)\|x\|_{p_1} \}
\backslash E ) \le  
A\exp(-B\|x\|).  
$$ 
\end{lemme} 

\begin{figure}[h!]
\begin{flushleft}
\input{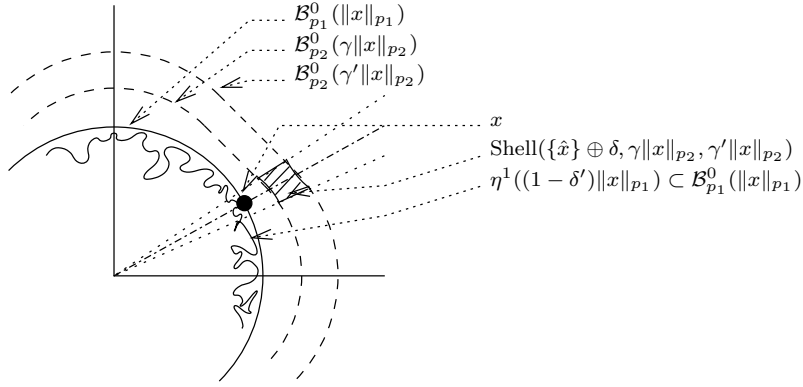}
\caption{Initialization of the spread.} 
\label{DESSinit}
\end{flushleft}
\end{figure}

\begin{proof} Let $\delta>0$ and choose any $0<\delta'<\delta$. 

We first need to introduce a certain number of parameters: Let
\begin{eqnarray}
0< & \rho & < \frac1{1-\delta'/2}-1 \label{rhoa} 
\end{eqnarray}
By Lemma~\ref{contunif}, we can then choose $\theta_1$ such that
\begin{equation}
\label{theta1a}
\|\hat x -\hat y\|_{p_2} \le \theta_1 
\quad \Longrightarrow  \quad (1-\rho)
\frac{\|x\|_{p_2}}{\|x\|_{p_1}} \le \frac{\|y\|_{p_2}}{\|y\|_{p_1}} \le (1 +\rho)
\frac{\|x\|_{p_2}}{\|x\|_{p_1}}.
\end{equation}
Choose now $\gamma'_0>1$ and $\theta_0>0$ small enough to fulfill the
three following conditions:  
\begin{eqnarray*}
2(\gamma'_0-1)+3 \theta_0 & < & \theta_1 \\
1-\theta_0 & > & (1-\delta'/2)(1+\rho) \\
\gamma'_0-1 +\theta_0 & < & \frac{\delta-\delta'}{C_{p_1,p_2}}.
\end{eqnarray*}
Note that the second condition is allowed by the choice (\ref{rhoa}) for $\rho$.
As these conditions are monotone, they are still fulfilled for any
$\gamma' \in(1,\gamma'_0)$ and any $\theta \in (0,\theta_0)$. Choose then such a
$\theta$ and such a $\gamma'$, and choose $\alpha>0$ small enough to have:  
\begin{eqnarray}
2(1+\alpha)(\gamma'-1+\theta)+ \theta & < & \theta_1 
\label{C1}\\
(1+\alpha)(1-\theta)-\alpha \gamma' & > & (1-\delta'/2)(1+\rho)  
\label{C2}\\
(1+\alpha)(\gamma'-1 +\theta) & < & \frac{\delta-\delta'}{C_{p_1,p_2}}.
\label{C3}
\end{eqnarray}
Note that these conditions are allowed by the three previous ones. Finally,
choose any $1<\gamma<\gamma'$. 

\vspace{0.2cm} 
\underline{Step 0:} 
Suppose that $t(x) \le (1-\delta)\|x\|_{p_1}$. 
This implies that there exists a $p_2$-open finite path from the 
source $s_2$ to $x$, and by the classical estimate~(\ref{amasfini}) on the radius of finite
open clusters in supercritical 
percolation, there exist two strictly positive constants $A_0$ and $B_0$ such 
that 
\begin{equation} 
\label{equationa0} 
\forall x \in \Zd \quad \P(\{t(x) \le (1-\delta)\|x\|_{p_1}\} \backslash 
E_0 )\le A_0 \exp(-B_0\|x\|). 
\end{equation} 

\underline{Step 1:}
In this step, we use the typical spread of
first-passage percolation with parameter $p_1$ in a amount of time of 
$(1-\delta') \|x\|_{p_1}$. Note
$$E_1'= \left\{
\eta^1( (1-\delta') \|x\|_{p_1} ) 
\subset \mathcal{B}_{p_1}^0 \left((1-\delta'/2) \|x\|_{p_1} 
\right)  
\right\} \subset E_1.$$ 
The large
deviations result associated to the shape 
theorem~(\ref{shapeGD}) ensures that there exist two strictly positive
constants $A_1$  
and $B_1$ such that  
\begin{eqnarray}
\label{equationa1} 
\P((E_1')^c) & \le &  A_1 \exp(-B_1 \|x\|).
\end{eqnarray}

\underline{Step 2:}
In this step, we control the spread of first-passage
percolation with parameter $p_2$. Let us first prove the geometrical fact:  
\begin{equation} 
\label{geoa1} 
\left( \bigcup_{y \in \shell(\{\hat x\}\oplus\theta,
    \gamma\|x\|_{p_2},\gamma'\|x\|_{p_2})}  
\mathcal{B}_{p_2}^y((1+\alpha)\|y-x\|_{p_2}) \right) \cap 
\mathcal{B}_{p_1}^0((1-\delta'/2)\|x\|_{p_1})=\varnothing 
\end{equation} 
Note that since $y-x=(\|y\|_{p_2}-\|x\|_{p_2})\hat{y}+\|x\|_{p_2} 
(\hat{y}-\hat{x})$, we have, for every $y\in \shell(\{\hat x\}\oplus\theta,
\gamma\|x\|_{p_2},\gamma'\|x\|_{p_2})$, 
\begin{equation} 
\label{lecone} 
\|y-x\|_{p_2}\le 
\|y\|_{p_2}-\|x\|_{p_2}+\theta \|x\|_{p_2}\le (\gamma'-1+\theta)\|x\|_{p_2}.
\end{equation}  
Then, if $y\in \shell(\{\hat x\}\oplus\theta,
\gamma\|x\|_{p_2},\gamma'\|x\|_{p_2})$ and $z \in 
\mathcal{B}_{p_2}^y ( 
(1+\alpha)\|y-x\|_{p_2})$, we obtain first: 
\begin{eqnarray}
\|z\|_{p_2}
& \ge & \|y\|_{p_2} - \|z-y\|_{p_2} \nonumber 
 \ge  \|y\|_{p_2} - (1+\alpha)\|y-x\|_{p_2} \nonumber \\
& \ge & \|y\|_{p_2} - (1+\alpha) \left( \|y\|_{p_2}-\|x\|_{p_2}+\theta
  \|x\|_{p_2} \right) \quad \text{with (\ref{lecone})} \nonumber \\ 
& \ge & (1+\alpha)(1-\theta) \|x\|_{p_2} -\alpha \|y\|_{p_2} \nonumber \\
& \ge & \left( (1+\alpha)(1-\theta) -\alpha \gamma' \right) \|x\|_{p_2},
\label{mod2z}
\end{eqnarray}
and then:
\begin{eqnarray*}
\|\hat z-\hat x \|_{p_2}
& \le & \|\hat z-\hat y \|_{p_2} + \|\hat y-\hat x \|_{p_2} \\
& \le & \frac{2\|z-y\|_{p_2}}{\|y\|_{p_2}} + \theta \quad \text{ with
Lemma~\ref{direction}}\\
& \le & \frac{2(1+\alpha)\|x-y\|_{p_2}}{\|y\|_{p_2}} + \theta 
 \le  \frac{2(1+\alpha)(\gamma'-1+\theta)}{\gamma} + \theta   \quad
\text{with (\ref{lecone})}  \\
& \le & 2(1+\alpha)(\gamma'-1+\theta) + \theta < \theta_1
\text{ with assumption (\ref{C1}).}
\end{eqnarray*}
Thus, by definition~(\ref{theta1a}) of $\theta_1$, we have:
\begin{eqnarray*}
\|z\|_{p_1} 
& \ge & \left( \frac1{1+\rho} \right) \frac{\|z\|_{p_2}}{\|x\|_{p_2}} \|x\|_{p_1}
\ge \frac{(1+\alpha)(1-\theta) -\alpha \gamma' }{1+\rho} \|x\|_{p_1}\quad
\text{ with (\ref{mod2z})} \\ 
& > & (1-\delta'/2) \|x\|_{p_1}\quad \text{ with assumption (\ref{C2})},
\end{eqnarray*}
which proves 
inclusion~(\ref{geoa1}).

Now, if we denote 
$$E_2'=\bigcap_{y \in 
  \shell(\{\hat x\}\oplus\theta, \gamma\|x\|_{p_2},\gamma'\|x\|_{p_2}) \cap
  C_{p_2}^\infty }    
\left\{ 
\begin{array}{c}
x \stackrel{p_2}{\communique} y \text{ in } 
\mathcal{B}_{p_2}^y((1+\alpha)\|y-x\|_{p_2}) \\
\text{by a path of length smaller} \\
\text{than } (1+\alpha)\|y-x\|_{p_2}
\end{array} 
\right\}
$$ 
Lemma~\ref{chemin-pas-loin} ensures that there exist two strictly positive 
constants $A_2,B_2$ such that  
\begin{equation} 
\label{equationa2} 
\forall x \in \Zd \quad \P \left( \{x \in  C_{p_2}^\infty\} \backslash E_2'
\right) 
\le C_d (\gamma'\|x\|_{p_2})^d A_2 \exp(-B_2 (\gamma-1)\|x\|).
\end{equation}

\underline{Step 3:} 
To conclude, it only remains to see that:
\begin{equation}
\label{geoa3}
\forall y \in \shell(\{\hat x\}\oplus\theta,
\gamma\|x\|_{p_2},\gamma'\|x\|_{p_2}) \quad 
(1+\alpha)\|y-x\|_{p_2} 
\le (\delta-\delta')\|x\|_{p_1}.
\end{equation}
Indeed, we have, for any $y \in \shell(\{\hat x\}\oplus\theta,
\gamma\|x\|_{p_2},\gamma'\|x\|_{p_2})$:  
\begin{eqnarray*}
(1+\alpha)\|y-x\|_{p_2} 
& \le & (1+\alpha)(\gamma'-1+\theta)\|x\|_{p_2} \quad \text{ with
(\ref{lecone})} \\
& \le & (1+\alpha)(\gamma'-1+\theta)C_{p_1,p_2}\|x\|_{p_1} \\
& \le & (\delta-\delta')\|x\|_{p_1}  \quad \text{ with assumption (\ref{C3})}.
\end{eqnarray*}
Now, if $E'=E_0 \cap
E'_1 \cap E'_2$, then 
Equation~(\ref{geoa1}) and 
inclusion~(\ref{geoa3}) imply that $\{t(x) \le (1-\delta)\|x\|_{p_1} \}\cap
E' \subset \{t(x) \le (1-\delta)\|x\|_{p_1} \}\cap E$, and thus
\begin{eqnarray*}
 \P( \{ t(x) \le (1-\delta)\|x\|_{p_1} \}  
\backslash E ) 
& \le & \P( \{ t(x) \le (1-\delta)\|x\|_{p_1} \} \backslash E_0) \\
&& + \P( \{ t(x) \le (1-\delta)\|x\|_{p_1} \} \backslash E'_1) \\
&& + \P( \{ t(x) \le (1-\delta)\|x\|_{p_1} \} \backslash E'_2 \cap E_0).
\end{eqnarray*}
Equations~(\ref{equationa0}), (\ref{equationa1}), 
(\ref{equationa2}) and the fact that if $t(x) \le (1-\delta)\|x\|_{p_1}$,
then $t(x)+(\delta-\delta')\|x\|_{p_1} \le (1-\delta')\|x\|_{p_1}$ give the
desired result.  
\end{proof}

\subsection{Typical progression of the stronger infection from one shell to the next one} 
$\;$

In this subsection, we forget for a moment the competition model, and study
the progression 
of one infection with parameter $p_2$. For simplicity, we omit, only in this subsection, the subscript $p_2$.
In the
next lemma, we want to bound the
minimal time needed 
for the infection to colonize the big $\shell(T,
(1+h)r, (1+h)^2 r)$ 
from the small $\shell(S, r, (1+h)r)$. 

\begin{lemme} 
\label{progressionter}
Let $\phi \in (0,2]$, $h \in (0,1)$ and $\alpha \in (1,2)$ be fixed parameters such that
\begin{equation}
(1+h)^2 (1+ \phi) - (1+h) 
< \alpha h <2  \label{C11}.
\end{equation}

For any $S$ and $T$ subsets of $\mathcal{S}$ and for any $r>0$, we define
the following event $E=E(S,T,r)$:
"Any point in the big $\shell(T, (1+h)r, (1+h)^2 r)\cap C^\infty$ is linked to a point in the small 
$\shell(S, r, (1+h)r)$ by an open path whose length is
less than $\alpha h r$.'' Two subsets $S$ and $T$ of $\mathcal{S}$ are said to be ``good'' if 
$$\forall \hat z \in T \quad \exists \hat v_z \in \mathcal{S}\text{ such that }
\{\hat v_z\} \oplus \frac{\phi}{2} \subset S \text{ and } \|\hat z -\hat v_z\|
\le \phi.$$

Then there exist two strictly positive constants $A$ and $B$, only depending
on $\phi, h, \alpha$, such that for any $r>0$ and 
for any two ``good'' subsets $S$ and $T$ of $\mathcal{S}$, 
we have $\P(E^c) \le A \exp(-Br)$.

Moreover, we can assume  that all the infection paths needed in $E$ are
completely  
included in the bigger $\shell(T\oplus(2\alpha h), [1- 3\phi] (1+h)r, \infty)$.
\end{lemme}

\begin{figure}[h!]
\begin{flushleft}
\input{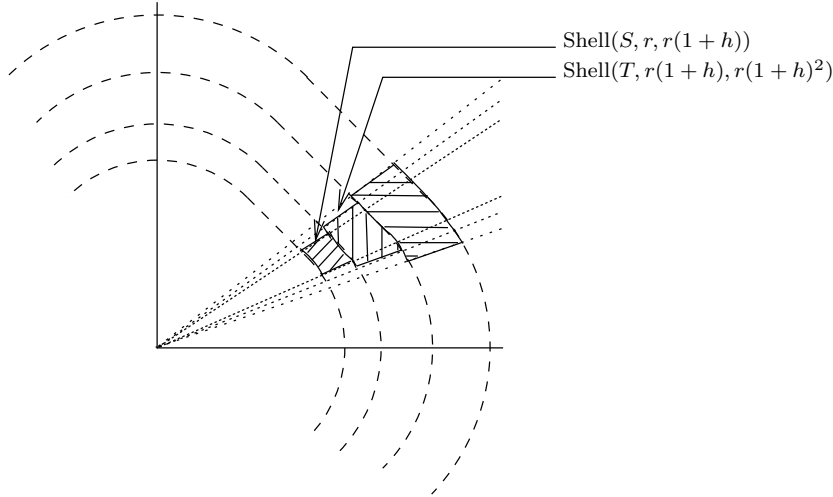}
\caption{Scheme of progression of a single infection from one shell to the other. A possible third larger shell is also drawn.} 
\label{DESSprog}
\end{flushleft}
\end{figure}

\begin{proof}
Let $\phi \in (0,2]$, $h \in (0,1)$ and $\alpha \in (1,2)$ be fixed parameters satisfying
Equation~(\ref{C11}) and choose, in this order, $\alpha'>1$,
$\epsilon>0$ and $\rho>0$ such that
\begin{eqnarray}
(1+h)^2 (1+ \phi) - (1+h-2\rho) 
& \le & \alpha' h < \alpha h, \label{C21} \\
\frac{2\rho \alpha' h }{1+h -2 \rho} 
& \le &  \frac{\phi}{2}, \label{C22} \\
h-2\rho-\rho \alpha' h & > & 0 ,\label{C24} \\
(1+\epsilon)(1+\rho)\alpha' & \le & \alpha.\label{C25}
\end{eqnarray}
Take any two ``good'' subsets $S$ and $T$ of $\mathcal{S}$. 
For any $z \in \shell(T, (1+h)r, (1+h)^2 r)$, we can choose $\hat v_z \in \mathcal{S}$ such that
\begin{equation}
\begin{array}{l}
\{\hat v_z\} \oplus \frac{\phi}{2} \subset S \quad \text{and} \quad \|\hat z -\hat v_z\|
\le \phi,
\end{array}
\nonumber
\end{equation}
and we set $v_z=(1+h-2\rho)r \hat v_z$. 
Let us first estimate $\|z-v_z\|$: on the one hand,
\begin{eqnarray}
\|z-v_z\| & \le & \|z-\|z\|\hat{v}_z\|+| \, \|z\|- (1+h-2\rho)r| \nonumber \\
& \le & \|z\| \phi + \|z\|- (1+h-2\rho)r \nonumber \\
& \le &  \|z\| (1+ \phi )-(1+h-2\rho)r \label{modulezvMaj} \\
& \le & [(1+h)^2 (1+ \phi) - (1+h-2\rho) ] r 
\le \alpha' h r  \text{ thanks to (\ref{C21})} \label{majax},
\end{eqnarray}
and, on the other hand,
\begin{eqnarray} 
\|z-v_z\| & \ge & \|z\|-\|v_z\| \ge 2 \rho r  \label{modulezvMin}.
\end{eqnarray}

\vspace{0.2cm}
\underline{Idea of the proof:} 
The idea of the proof is the following. Take $z$ in $C^\infty_{p_2}$ and in $\shell(T, (1+h)r,
(1+h)^2 r)$.
The ball $\mathcal{B}^{v_z}( \rho \|z-v_z\| )$ is included in the $\shell(S, r, (1+h)r)$ and in 
the ball $\mathcal{B}^{z}( (1+\rho) \|z-v_z\| )$. As it is of radius of order $r$, it contains with high probability some point of the infinite cluster, and this point should be with high probability, thanks to Lemma~\ref{chemin-pas-loin}, linked to $z$ by an open path inside
$\mathcal{B}^{z}( (1+\epsilon)(1+\rho) \|z-v_z\| )$ with length less
than $(1+\epsilon)(1+\rho) \|z-v_z\| )$. We
chose the parameters to ensure that $(1+\epsilon)(1+\rho) \|z-v_z\| \le
\alpha h r$. 
It will then only remain to control the positions of the points in the union of
the
$\mathcal{B}^{z}( (1+\epsilon)(1+\rho)  \|z-v_z\|)$. Let us
make all this more precise.

\vspace{0.2cm}
\underline{Geometrical facts:} 
Let us first note that, by the triangular inequality,  
\begin{equation}
\label{fact1}
\forall z \in \shell(T, (1+h)r, (1+h)^2 r) \quad \mathcal{B}^{v_z}( \rho \|z-v_z\|) \subset \mathcal{B}^{z}( (1+\rho)  \|z-v_z\|).
\end{equation}
Let us see now that 
\begin{equation}
\label{fact2}
\forall z \in \shell(T, (1+h)r, (1+h)^2 r) \quad \mathcal{B}^{v_z}( \rho
\|z-v_z\|) \subset \shell(S , r, (1+h)r). 
\end{equation}
Let $u \in \mathcal{B}^{v_z}( \rho \|z-v_z\| )$, then, by Lemma~\ref{direction},
\begin{eqnarray*} 
\|\hat u-\hat v_z\| & \le & \frac{2\rho \|z-v_z\|}{\|v_z\|}  
\le \frac{2\rho \alpha' h}{1+h-2\rho} \text{ by Equation~(\ref{majax}) and
  definition of } v_z \\
& \le & \frac{\phi}{2} \text{ thanks to Equation~(\ref{C22})}
\end{eqnarray*}
and thus $\hat u \in S$.
For the norm of $u$, by definition of $v_z$ and Equation~(\ref{majax}), we have:
\begin{eqnarray*}
\|v_z\|- \rho \|z-v_z\| \le & \|u\| & \le  \|v_z\|+
\rho \|z-v_z\| \\ 
(1+h-2\rho)r - \rho \alpha' h r \le &  \|u\| & \le  (1+h-2\rho)r +
\rho \alpha'  h r  \\
r \le & \|u\| & \le (1+h)r, 
\end{eqnarray*}
thanks to Equations~(\ref{C24}) and (\ref{C21}).
This proves the second geometrical fact~(\ref{fact2}).

\vspace{0.2cm}
\underline{Probabilistic estimates:} 
We can then estimate the probability of $E$.
Note first
$$E_1=\bigcup_{z \in \shell(T,(1+h) r, (1+h)^2r)} \left\{
\mathcal{B}^{v_z} \left(\rho\|z-v_z\|\right) \cap C^\infty_{p_2} =\varnothing
\right\}.$$
By estimate (\ref{modulezvMin}), we know that $\|z-v_z\|\ge 2\rho r$; moreover, for each 
$z\in\shell(T,(1+h) r, (1+h)^2r)$, the point $v_z$ is in $ \shell(S, r, (1+h)r)$. Thus, using the estimate on the holes of the infinite cluster~(\ref{amasinfini}), there exist two strictly positive constants $A_1$ and $B_1$ such that
for every ``good'' $S$ and $T$, for every $r>0$,
\begin{eqnarray*}
\P(E_1) & \le & 
 \P \left(
\bigcup_{v \in \shell(S, r, (1+h)r)} \left\{
\mathcal{B}^{v} \left( 2\rho^2 r \right) \cap C^\infty_{p_2} =\varnothing
\right\}
\right) \\
& \le & C_d [(1+h)r]^d A_1 \exp (-B_1  2\rho^2 r).
\end{eqnarray*}
Then, note
$$E_2=
\bigcup_{ 
\begin{subarray}{c}
{z \in \shell(T, (1+h)r, (1+h)^2 r) }\\
u \in \mathcal{B}^{v_z}( \rho \|z-v_z\|)
\end{subarray}}
\left\{
\begin{array}{c}
u \in C^\infty_{p_2}, \; 
z \in C^\infty_{p_2}, \text{ and $u$ is not linked to }  \\
z \text{ by an open path of length smaller} \\
\text{than } (1+\epsilon)(1+\rho)  \|z-v_z\|\\
\text{inside } \mathcal{B}^{z}((1+\epsilon)(1+\rho)  \|z-v_z\|) 
\end{array}
\right\}.
$$
By Lemma~\ref{chemin-pas-loin}, 
Equations~(\ref{modulezvMin}) and (\ref{majax}), there exist two strictly positive 
constants $A_2$ and $B_2$ such that 
for every ``good'' $S$ and $T$, for every $r>0$,
\begin{eqnarray*}
\P(E_2)& \le & 
\sum_{ 
\begin{subarray}{c}
{z \in \shell(T, (1+h)r, (1+h)^2 r) }\\
u \in \mathcal{B}^{v_z}( \rho \|z-v_z\|)
\end{subarray}
}
\P
\left(
H(u,z,(1+\epsilon)(1+\rho)-1)^c
\right)
 \\
& \le & \sum_{ 
\begin{subarray}{c}
{z \in \shell(T, (1+h)r, (1+h)^2 r) }\\
u \in \mathcal{B}^{v_z}( \rho \|z-v_z\|)
\end{subarray}
} A_2 \exp (-B_2  (1+\epsilon)(1+\rho) \|z-v_z\|)\\
& \le & C_d [(1+h)^2r]^d \times C_d (\rho \alpha' h r)^d \times A_2 \exp (-B_2 (1+\rho)2 \rho r).
\end{eqnarray*}

\vspace{0.2cm}
\underline{Conclusion:}
For every $z \in \shell(T, (1+h)r, (1+h)^2 r)$, thanks to (\ref{majax})
and (\ref{C25}), one has $(1+\epsilon)(1+\rho) 
\|z-v_z\|_{p_2}\le \alpha h r $. This, combined with geometrical facts~(\ref{fact1}) and (\ref{fact2}), implies that $E^c \subset E_1 \cup E_2$, which proves the exponential estimate of the lemma. 

\vspace{0.2cm}
\underline{Control of the infection paths:} 
It remains to estimate the minimal room needed to perform this infection, or in
other words to control
$$\bigcup_{z \in \shell(T, (1+h)r, (1+h)^2 r)} \mathcal{B}^{z}(
(1+\epsilon)(1+\rho) \|z-v_z\| ).$$ 
Let $z \in \shell(T , (1+h)r, (1+h)^2 r)$ and $u \in \mathcal{B}^{z}(
(1+\epsilon)(1+\rho) \|z-v_z\|)$. We have:
\begin{eqnarray*}
\|u\| & \ge & \|z\| - (1+\epsilon)(1+\rho) \|z-v_z\|
\\
& \ge & (1+\epsilon)(1+\rho)(1+h-2\rho)r -
[(1+\epsilon)(1+\rho)(1+\phi)-1]\|z\| \\
& & \quad \quad \text{
thanks to (\ref{modulezvMaj})} \\
& \ge & (1+\epsilon)(1+\rho)(1+h-2\rho)r -
[(1+\epsilon)(1+\rho)(1+\phi)-1](1+h)^2r \\
& \ge & [1-3\phi] (1+h)r.
\end{eqnarray*}
The last inequality is obtained by looking at the the limit of the right-hand side 
term, when $\epsilon$ and $\rho$ tend to $0$, and  by decreasing if
necessary $\epsilon$ and $\rho$. Finally, by applying Lemma~\ref{direction} and then Inequality~(\ref{majax}), we have
$$
\|\hat u -\hat z\|
\le \frac{2\|u-z\|}{\|z\|} 
\le  \frac{2 (1+\epsilon)(1+\rho)\alpha' h}{(1+h)^2}
 \le 2\alpha h.
$$
Thus $u \in \shell(T\oplus(2\alpha h), r_{min},\infty)$, which ends the proof of the lemma.
\end{proof}

\subsection{Final step: proof of Theorem~\ref{avantage2}}
We come back now to the competition context, with a weaker infection with parameter $p_1$ and a stronger infection with parameter $p_2>p_1$.

\begin{proof}
Let $\delta>0$. 

\vspace{0.2cm}
\underline{Idea of the proof:} 
The idea is quite natural: start the progression by the initialization Lemma~\ref{petit-bout-de-coquille}, and apply recursively the progression Lemma~\ref{progressionter} until the stronger infection surrounds the weaker one. The point is to ensure that this progression is not disturbed by the spread of the weaker infection.

\vspace{0.2cm}
\underline{Step 0. Choice of constants:}
Remember that $C_{p_1,p_2}<1$ and choose:
\begin{eqnarray}
\delta'>0 & \text{such that} & \delta'<\delta \text{ and } 
\delta'<1-C_{p_1,p_2}, \label{D1} 
\\
\rho>0 & \text{such that} & (1+\rho)(1-\delta')<1. \label{D6}
\end{eqnarray}
By Lemma~\ref{contunif}, there exists ${\theta}>0$ such that for any $x,y \in \Zd \backslash \{0\}$, we have:
\begin{equation}
 \|\hat x -\hat y\| \le {\theta}  \Rightarrow (1-\rho)
\frac{\|x\|_{p_2}}{\|x\|_{p_1}} \le \frac{\|y\|_{p_2}}{\|y\|_{p_1}} \le (1
+\rho) \frac{\|x\|_{p_2}}{\|x\|_{p_1}}. \label{D9}
\end{equation}
Choose then $h$  and $\alpha$ such that:
\begin{eqnarray}
0<h<1 & \text{such that} & (1+h)C_{p_1,p_2}<1-\delta', 
\nonumber \\
1<\alpha< 2 & \text{such that} & \alpha > 1+h  \text{ and }  \alpha
C_{p_1,p_2}<1-\delta', \label{D3} \\
2\alpha h < {\theta}. & \label{D2}
\end{eqnarray}
The first condition is allowed by condition~(\ref{D1}) on $\delta'$, and
allows itself the choice~(\ref{D3}) for $\alpha$. We obtain~(\ref{D2}) by
decreasing $h$ if necessary.
Let $\gamma'_0>1$ and $\theta_0>0$ be given by
Lemma~\ref{petit-bout-de-coquille}.
Choose $\gamma',\gamma,\epsilon$ and $\phi$ in the
following manner:
\begin{eqnarray}
1<\gamma'<\gamma'_0 & \text{such that} & \alpha \gamma'
C_{p_1,p_2}<1-\delta', \label{D4} \\
1<\gamma<\gamma' & \text{such that} & \gamma=\frac{\gamma'}{1+h}, \label{D5}
\\
\epsilon>0 & \text{such that} & 
\left\{
\begin{array}{ll}
& \alpha C_{p_1,p_2}(1+\epsilon)<1, \\
\text{and} & (1+\epsilon)(1+\rho)(1-\delta')<\gamma,  
\end{array} 
\right.
\label{D7} \\
0<\phi<\theta_0 & \text{such that} & 
\left\{
\begin{array}{ll}
& \alpha C_{p_1,p_2}(1+\epsilon)<1 -3 \phi, \\
\text{and} & (1+\epsilon)(1+\rho)(1-\delta')<\gamma (1-3\phi), \\ 
\text{and} & (1+h)^2 (1+ \phi) - (1+h) < \alpha h. 
\end{array} 
\right.
\label{D8}  
\end{eqnarray}
Note that condition~(\ref{D4}) is allowed by the choice~(\ref{D3}), and
condition~(\ref{D5}) is obtained by decreasing $h$ if
necessary. Conditions~(\ref{D7}) are respectively permitted by (\ref{D3})
and (\ref{D6}), and allow the first two conditions on $\phi$. The last
condition in (\ref{D8}) is allowed by (\ref{D3}) and (\ref{D2}).
Choose now $K\ge2$ large enough to have for every $k \ge K$
\begin{equation}
\label{D10}
C_{p_1,p_2}(1+\epsilon) \left[ (1-\delta') {C_{p_2,p_1}}
  + \alpha \gamma ((1+h)^{k-1}-1) \right]< \gamma [1-3\phi](1+h)^{k-1},
\end{equation}
which is allowed by (\ref{D8}). By decreasing $\phi$ if necessary, we can assume, thanks to~(\ref{D2}), that 
\begin{equation}
\label{D12}
(1+K)\frac{\phi}{2}+ 2 \alpha h < {\theta}.
\end{equation}

\vspace{0.2cm}
\underline{Step 1. Initialization of the spread:}
By Lemma~\ref{petit-bout-de-coquille}, there exist two strictly positive
constants $A_1$ and $B_1$ such that for every $x \in \Zd \backslash
\{0\}$, we have
\begin{equation}
\P( \{ t(x) \le (1-\delta)\|x\|_{p_1} \}
\backslash \{ E_1(x) \cap \{x\in C_{p_2}^\infty\} \} ) \le  
A_1\exp(-B_1\|x\|), \label{D11}
\end{equation}
where we use the following notations: 
\begin{eqnarray*}
E_1^1(x) & = & \{\eta^1( (1-\delta')\|x\|_{p_1}) \subset
\mathcal{B}_{p_1}^0(\|x\|_{p_1}) \}, \\
E_1^2(x) & = & \left\{\eta^2((1-\delta')\|x\|_{p_1}) \supset 
\left( \shell\left(\{\hat x\}\oplus\frac{\phi}{2}, \gamma\|x\|_{p_2},\gamma'\|x\|_{p_2}\right)
  \cap C_{p_2}^\infty \right) \right\}, \\  
E_1(x) & = & E_1^1(x) \cap E_1^2(x).
\end{eqnarray*}
Thus, if $t(x) \le (1-\delta)\|x\|_{p_1}$, then at the slightly larger
time $t_1(x)=(1-\delta') \|x\|_{p_1}$, the first shell 
$$S_1(x)=C_{p_2}^\infty \cap \shell\left( \{\hat x\}\oplus\frac{\phi}{2},
\gamma\|x\|_{p_2},\gamma'\|x\|_{p_2}\right) $$
is with high probability colonized by the $p_2$-infection. 

We want now to extend this colonization to larger and larger shells by
applying  recursively Lemma~\ref{progressionter}.

\vspace{0.2cm}
\underline{Notations:}
We still need to introduce a certain number of notations, inspired by Lemma~\ref{progressionter}:
$$
\begin{array}{|l|}
\hline
k=1 \\
\hline
r_1=\gamma \\
A_1(x)=\{\hat x\} \oplus \frac{\phi}{2} \\
S_1(x)=C^\infty_{p_2} \cap \shell(A_1(x), \gamma \|x\|_{p_2}, \gamma'\|x\|_{p_2}) \\
t_1(x)=(1-\delta') \|x\|_{p_1} \\
\hline
k \ge 2 \\
\hline
r_k=(1+h)^{k-1}r_1 
\quad \text{and} \quad r_k^{min}=[1-3\phi](1+h)r_{k-1} \\
A_k(x)=A_{k-1}(x)\oplus \frac{\phi}{2} \\
S_k(x)=
C^\infty_{p_2} \cap \shell(A_k(x), r_k\|x\|_{p_2},r_{k+1}\|x\|_{p_2}) \\
t_k(x)=t_{k-1}(x)+h \alpha \gamma r_{k-1}\|x\|_{p_2} =(1-\delta') \|x\|_{p_1} +
\alpha[(1+h)^{k-1}-1] \|x\|_{p_2} \\
\quad \quad =t_1(x)+ \alpha(r_k-r_1) \|x\|_{p_2}\\
\hline
\end{array}
$$
Define also the following events, for $k \ge 2$ and $x \in \Zd \backslash
\{0\}$:
\begin{eqnarray*}
E_k^1 (x) & = & \{\eta^1(t_k(x)) \subset
\mathcal{B}_{p_1}^0((1+\epsilon)t_k(x)) \}, \\ 
E_k^2 (x)& = & \{\eta^2(t_k(x)) \supset S_k(x) 
\}, \\ 
E_k (x)& = & E_k^1 (x) \cap E_k^2(x).
\end{eqnarray*}
The aim is the following: we want to apply Lemma~\ref{progressionter} to prove that if $E_k^2 (x)$ is
realized, then with high probability $E_{k+1}^2 (x)$ is also realized. But we
need first to control the spread of the slow $p_1$-infection, and to see
that it will not disturb the spread of the fast $p_2$-infection from
$S_k(x)$ to $S_{k+1}(x)$.

\vspace{0.2cm}
\underline{Step 2. Rough control of the slow $p_1$-infection:} Let us prove that
there
exist two strictly positive constants $A_2$ and $B_2$ such that
\begin{equation}
\label{EE2}
\forall x \in \Zd \backslash\{0\}  \quad 
\P \left( \complement \left( \miniop{}{\bigcap}{k \ge 2} E_k^1(x) \right) \right) \le A_2 \exp (-B_2 \|x\|).
\end{equation} 
Indeed, by the large deviation result~(\ref{shapeGD}), for any $x \in \Zd \backslash\{0\}$, we have:
\begin{eqnarray*}
&& \P \left( \complement \left( \bigcap_{k \ge 2} E_k^1(x) \right) \right) \\
& \le & \sum_{k \ge 2} \P \left( \eta^1(t_k(x)) \not\subset
  \mathcal{B}_{p_1}((1+\epsilon)t_k(x))) \right)  
 \le  \sum_{k \ge 2} A \exp \left( -B t_k(x) \right) \\
& \le & \sum_{k \ge 2} A \exp \left( -B 
  \left[(1-\delta')\|x\|_{p_1}+\alpha((1+h)^{k-1}-1)\gamma
    \|x\|_{p_2}\right] \right) \\ 
& \le & A \exp \left( -B \left[ (1-\delta')\|x\|_{p_1} -\alpha \gamma
\|x\|_{p_2} \right] \right) \sum_{k \ge 2} \exp \left( -B
\alpha(1+h)^{k-1}\gamma \|x\|_{p_2} \right).
\end{eqnarray*}

1. As there exists $B'>0$ such that
$\forall k \ge 2, \;B \alpha(1+h)^{k-1}\gamma \ge B' k,$
the last sum is bounded by
\begin{eqnarray*}
\sum_{k \ge 2} \exp \left( -B
\alpha(1+h)^{k-1}\gamma \|x\|_{p_2} \right)
& \le & \sum_{k \ge 2} \exp \left( -B' k\|x\|_{p_2} \right) \\
& \le & \frac{\exp \left( -2B' \|x\|_{p_2} \right)}{1-\exp \left( -B'
    \|x\|_{p_2} \right)} \le A' \exp \left( -2B' \|x\|_{p_2} \right)
\end{eqnarray*}
with $A'>0$ because $\displaystyle \inf \{\|x  \|_{p_2}: \; x \in \Zd \backslash\{0\} \} >0$.

2. For the first factor, we have
\begin{eqnarray*}
(1-\delta')\|x\|_{p_1} -\alpha \gamma \|x\|_{p_2}
& = & (1-\delta'-\alpha \gamma C_{p_1,p_2})\|x\|_{p_1} + \alpha \gamma
(C_{p_1,p_2} \|x\|_{p_1}-\|x\|_{p_2}) \\
& \ge & (1-\delta'-\alpha \gamma C_{p_1,p_2})\|x\|_{p_1} \text{ by
  definition of } C_{p_1,p_2} \\
& \ge & B''\|x\|_{p_1},
\end{eqnarray*}
with $B''>0$ thanks to conditions  (\ref{D4}) and (\ref{D5}). This proves
(\ref{EE2}).  

But we will also need a more precise control of this slow infection in
order to prevent it from bothering the fast one while applying
Lemma~\ref{progressionter}. 

\vspace{0.2cm}
\underline{Step 3. More precise control of the slow $p_1$-infection for large times:}
Remember that $K$ was defined in (\ref{D10}). Let us prove the following
geometrical fact:
\begin{equation}
\label{EE3}
\forall k \ge K \quad \forall x \in \Zd \backslash\{0\} \quad
\mathcal{B}_{p_1}^0((1+\epsilon)t_k(x)) 
\subset\mathcal{B}_{p_2}^0(r_k^{min}\|x\|_{p_2}).
\end{equation}
Let $k \ge K $, $x \in \Zd \backslash\{0\}$ and $u \in \mathcal{B}_{p_1}^0((1+\epsilon)t_k(x))$. Then:
\begin{eqnarray*}
\|u\|_{p_2} 
& \le & C_{p_1,p_2} \|u\|_{p_1} \le C_{p_1,p_2} (1+\epsilon)t_k(x) \\
& \le & C_{p_1,p_2}
(1+\epsilon)[(1-\delta')\|x\|_{p_1}+\alpha \gamma ((1+h)^{k-1}-1)
\|x\|_{p_2}] \\
& \le & C_{p_1,p_2} (1+\epsilon) \left[ (1-\delta') {C_{p_2,p_1}}
  + \alpha \gamma ((1+h)^{k-1}-1)  \right] \|x\|_{p_2} \\
& \le & [1-3\phi](1+h)^{k-1}\gamma\|x\|_{p_2}=r_k^{min}\|x\|_{p_2},
\end{eqnarray*}
thanks to (\ref{D10}), which proves (\ref{EE3}).

\vspace{0.2cm}
\underline{Step 4. More precise control of the slow $p_1$-infection in the early stage of the process:} \\
To look at the $p_1$-infection in the early stage of the process, we need
to focus on a small 
cone around $\hat x$ in order to control more precisely the discrepancy
between the two norms $\|.\|_{p_1}$ and $\|.\|_{p_2}$. 
Let us see that for every $k$, for every $x\in \Zd \backslash \{0\}$ and for every $z \in
\Zd \backslash \{0\}$ 
\begin{equation}
(\|\hat z -\hat x \|_{p_2} \le {\theta} \text{ and } \|z\|_{p_1} \le
(1+\epsilon)t_k(x))  \Longrightarrow \|z\|_{p_2} \le
r_k^{min}\|x\|_{p_2}. \nonumber
\end{equation}
We recall that ${\theta}$ was defined in (\ref{D9}). Then, 
\begin{eqnarray*}
&& r_k^{min}\|x\|_{p_2}-\|z\|_{p_2} 
 \ge  r_k^{min}\|x\|_{p_2}-(1+\rho)\|z\|_{p_1}
\frac{\|x\|_{p_2}}{\|x\|_{p_1}}  \\
& \ge & (1-3\phi)r_k\|x\|_{p_2}- 
(1+\rho)(1+\epsilon)\frac{\|x\|_{p_2}}{\|x\|_{p_1}} t_k(x) \\ 
& \ge & (1-3\phi)r_k\|x\|_{p_2}-
(1+\rho)(1+\epsilon)\frac{\|x\|_{p_2}}{\|x\|_{p_1}}
(t_1(x)+\alpha(r_k-r_1)\|x\|_{p_2}). 
\end{eqnarray*}
As $(r_k)_{k}$ is increasing, the worst case is for $k=1$:
\begin{eqnarray*}
 r_k^{min}\|x\|_{p_2}-\|z\|_{p_2} 
& \ge & (1-3\phi)r_1\|x\|_{p_2} - (1+\rho)(1+\epsilon) \frac{\|x\|_{p_2}}{\|x\|_{p_1}} (1-\delta')\|x\|_{p_1} \\
& \ge & ((1-3\phi)\gamma - (1+\rho)(1+\epsilon)(1-\delta'))\|x\|_{p_2}>0
\end{eqnarray*}
thanks to Conditions~(\ref{D8}). Thus, thanks to Equation~(\ref{D12}), we
obtain that for every $k \le K $, for every $x\in \Zd \backslash \{0\}$ and for every $z \in \Zd \backslash \{0\}$
\begin{equation}
\label{EE4}
(\; \hat z \in A_k \oplus(2 \alpha h)  \quad \text{and} \quad \|z\|_{p_1} \le 
(1+\epsilon)t_k(x) \; )  \quad \Longrightarrow \quad \|z\|_{p_2} \le
r_k^{min}\|x\|_{p_2}. 
\end{equation}
 
\vspace{0.2cm}
\underline{Step 5. Control of the fast $p_2$-infection:}
Equations~(\ref{EE3}) and (\ref{EE4}) ensure that for every $k \ge 2 $, for
every $x\in \Zd \backslash \{0\}$, we have 
\begin{equation}
\label{EE6}
\mathcal{B}_{p_1}^0((1+\epsilon)t_k(x)) \; \cap \;
\shell( A_k(x) \oplus(2 \alpha h), r_k^{min} \|x\|_{p_2}, \infty)
= \varnothing. 
\end{equation}
Thus, the spread of the (single) fast $p_2$-infection from $S_{k-1}(x)$ to
$S_{k}(x)$, ensured 
by Lemma~\ref{progressionter}, is not disturbed by the slow $p_1$-infection
on the event $E_{k}^1(x)\cap E_{k-1}^2(x)$. Let $A_3$ and $B_3$ be the two
strictly positive constants given by Lemma~\ref{progressionter}; we apply the lemma with
\begin{eqnarray*}
S & = & A_{k-1}(x), \\
T & = & A_{k}(x)=A_{k-1}(x)\oplus \frac{\phi}{2} , \\
r & = & r_{k-1} \|x\|_{p_2}. 
\end{eqnarray*}
But we must first be sure that $S$ and $T$ are ``good''  subsets of $\mathcal S$, in the sense
$$ \forall z \in T \quad \exists v \in S \text{ such that } v \oplus \frac{\phi}{2} \subset S \text{ and } \|z-v\|_{p_2} \le \phi.$$
Indeed, let $k \ge 2$ and $z \in A_{k}(x)=A_{k-1}(x)\oplus \frac{\phi}{2}$: by definition, there exist
$w \in A_{k-1}(x)$ and $u_1 \in \mathcal{B}^0_{p_2}(\phi/2)$ such that $z=w+u_1$. But 
$A_{k-1}(x)=A_{k-2}(x)\oplus \frac{\phi}{2}$, where, for $k=2$, we set $A_0(x)=\{ \hat x\}$. So there exist
$v \in A_{k-2}(x)$ and $u_2 \in \mathcal{B}^0_{p_2}(\phi/2)$ such that $w=v+u_2$. Now, $z=v+u_1+u_2$ and 
\begin{itemize}
\item as $v \in A_{k-2}(x)\subset S$ and $S=A_{k-1}(x)=A_{k-2}(x)\oplus \frac{\phi}{2}$, we have  $v \oplus \frac{\phi}{2} \subset  S$,
\item as $u_1 \in \mathcal{B}^0_{p_2}(\phi/2)$ and $u_2 \in \mathcal{B}^0_{p_2}(\phi/2)$, we have $\|z-v\|_{p_2} =\|u_1+u_2\|_{p_2}\le \|u_1\|_{p_2}+\|u_2\|_{p_2} \le \phi$.
\end{itemize}

Thus any point in $S_{k}(x)$ can be infected by the $p_2$-infection from a
point in $S_{k-1}(x)$ in a time less than $\alpha h
r_{k-1}\|x\|_{p_2}=t_{k}(x)-t_{k-1}(x)$ using only paths inside $\shell( A_k(x) \oplus(2 \alpha h), r_k^{min} \|x\|_{p_2}, \infty)$, if it is not bothered by the slow $p_1$-infection. But on the event $E_{k}^1(x)$, this is ensured by Equation~(\ref{EE6}).
Thus, the application of  Lemma~\ref{progressionter} implies that
for 
any $x \in \Zd \backslash
\{0\}$, for every $k\ge2$,
\begin{equation}
\P((E_k^2(x))^c \backslash (E_k^1(x) \cap E_{k-1}^2(x))\le A_3 \exp(-B_3r_{k-1}\|x\|_{p_2}).
\nonumber
\end{equation}
Thus, 
\begin{eqnarray}
\sum_{k \ge 2} \P \left( (E_k^2(x))^c\cap (E_k^1(x) \cap
  E_{k-1}^2(x)\right) 
& \le & A_3 \sum_{k \ge 2} \exp(-B_3r_{k-1} \|x\|_{p_2}) \nonumber \\
& \le & A_3 \sum_{k \ge 2} \exp \left( -B_3(1+h)^{k-2}\gamma\|x\|_{p_2}
\right)\nonumber \\  
&  \le & A_4 \exp(-B_4\|x\|)\label{EE5} 
\end{eqnarray}
where $A_4$ and $B_4$ are two strictly positive constants.

\vspace{0.2cm}
\underline{Conclusion:} 
For $k$ large enough, the set $S_k(x)$ disconnects
$0$ from 
infinity, and thus the event $\bigcap_{k\ge1} E_k$ implies that the slow
$p_1$-infection is surrounded by the fast $p_2$-infection and thus
dies out. So, using (\ref{D11}), (\ref{EE2}) and (\ref{EE5}), we obtain:
\begin{eqnarray*}
&& \P(\mathcal{G}^1\cap \{t(x) \le (1-\delta) 
\|x\|_{p_1}\}) \\
& \le & \P \left( \{t(x) \le (1-\delta)
\|x\|_{p_1} \} \cap \bigcup_{k\ge 1} E_k(x)^c \right) \\
& \le & \P \left( \{t(x) \le (1-\delta)\|x\|_{p_1} \} \cap
(E_1(x) \cap\{x \in C_{p_2}^\infty\}) \right)
+ \P \left( \bigcup_{k\ge 1} (E_k^1(x))^c \right) \\
&& + \sum_{k \ge 2} \P \left( (E_k^2(x))^c\cap (E_k^1(x) \cap
  E_{k-1}^2(x)\right) \\
& \le & A\exp(-B\|x\|),
\end{eqnarray*}
which completes the proof.
\end{proof}

\section{Proof of the main Theorem~\ref{THEtheorem}}
\label{Section_proof_theo}

In all this section, $s_1$ and $s_2$ are two distinct sites in $\Zd$ and 
$\xi$ is the element of $S^{\Zd}$ where all sites are empty, but $\xi_{s_1}=\yellow$ and $\xi_{s_2}=\blue$. This initial configuration is now fixed.
We will thus, in the following, omit the explicit dependence in $\xi$.

Suppose that $0\le p \le q \le 1$.
In our competition process, the survival of the weaker -- \resp stronger -- infection is represented by the event $\mathcal{G}^1_{p,q}$ -- \resp $\mathcal{G}^2_{p,q}$  -- where, for $i=1,2$:
$$\mathcal{G}^i_{p,q}=\left\{ \sup_{t \ge 0} |\eta^i_{p,q}(t)|=+\infty\right\}.$$
The main Theorem~\ref{THEtheorem} can be reformulated now in the following form:
\begin{theorem}
\label{THETHEtheorem}
Let $ q>p_c$. The set of parameters $p$ such that $p<\min(q,\pcfleche)$ and
$\P(\mathcal{G}^1_{p, q}\cap \mathcal{G}^2_{p, q})>0$ is at most denumerable.
\end{theorem}
The corresponding Conjecture~\ref{Introconj} can be formulated as follows:
\begin{conjecture} 
Let $ q>p_c$  and $p<\min(q,\pcfleche)$. Then 
$\P \left( \mathcal{G}^1_{p,q} \cap \mathcal{G}^2_{p,q} \right) =0$.
\end{conjecture}

\subsubsection*{Proof of Theorem~\ref{THETHEtheorem}} 
It strongly relies on Propositions~\ref{Introstrictecomp} and~\ref{Introslowspeed} and the coupling arguments that are also used are widely inspired by
the proof of Häggström and Pemantle~\cite{Haggstrom-Pemantle-2}.

\vspace{0.2cm}
\underline{Step 1.} 
Let us prove that if $p<q<\min(r,\pcfleche)$, then
$\P(\mathcal{G}^1_{p,r}\cap \mathcal{G}^2_{q,r})=0$.

\vspace{0.2cm}
Since, by the coupling Lemma~\ref{yal}, $\mathcal{G}^2_{p,r}\subset \mathcal{G}^2_{q,r}$, we have
$\mathcal{G}^1_{p,r}\cap \mathcal{G}^2_{q,r}=(\mathcal{G}^1_{p,r}\cap \mathcal{G}^2_{p,r})\cap \mathcal{G}^2_{q,r}$.
So, we can assume that $\mathcal{G}^1_{p,r}\cap \mathcal{G}^2_{p,r}$ occurs and prove that $\mathcal{G}^2_{q,r}$ can not happen.
By Proposition~\ref{Introslowspeed}, we have 
$$
\miniop{}{\limsup}{t\to +\infty}\frac{|\eta^2_{p,r}(t)|_{p}}t\le 1, 
\text{ which implies }
\miniop{}{\limsup}{t\to +\infty}\frac{|\eta^2_{q,r}(t)|_{p}}t\le 1.
$$
Indeed, by the coupling Lemma~\ref{yal}, $\eta^2_{q,r}(t) \subset \eta^2_{p,r}(t)$.
Now, by Proposition~\ref{Introstrictecomp}, we have
$$\miniop{}{\limsup}{t\to +\infty}\frac{|\eta^2_{q,r}(t)|_{q}}t\le C_{p,q}.$$
On the other hand, $\mathcal{G}^1_{q,r}\subset\{s_1\in C^\infty_{q}\}$, so using Lemma~\ref{comparaison} and Lemma~\ref{shape} together, we get
$$\miniop{}{\liminf}{t\to +\infty}\frac{|\eta^1_{q,r}(t)\cup \eta^2_{q,r}(t)|_{*,q}}t\ge 1.$$
Now, let $t$ be large enough to ensure that 
$$\frac{|\eta^2_{q,r}(t)|_{p}}t\le\frac{C_{p,q}+2}3=\alpha \;
\text{ and } \;
\frac{|\eta^1_{q,r}(t)\cup \eta^2_{q,r}(t)|_{*,q}}t\ge \frac{2C_{p,q}+1}3=\beta.$$
Then every point $x$ such that $x\in C^\infty_{q}$ and $\alpha <\|x\|_{q}< \beta$
belongs to $\eta^1_{q,r}(t)\backslash\eta^2_{q,r}(t)$, which prevents the occurrence of the event
$\mathcal{G}^2_{q,r}$.

\vspace{0.2cm}
\underline{Step 2.} 
Let $q>0$. Let us prove that $\P$ almost surely,
there exists at most one value $p \le \min(q,\pcfleche)$ such that 
  $\mathcal{G}^1_{p,q} \cap \mathcal{G}^2_{p,q}$ occurs.
\vspace{0.2cm}

Assume that there exist $p$ and $p'$ with
$p<p'\le  q$ and such that $\mathcal{G}^1_{p, q}\cap \mathcal{G}^2_{p, q}$
and $\mathcal{G}^1_{p', q}\cap \mathcal{G}^2_{p', q}$ are satisfied. Denote by $A$ this event.
Let $r$ and $s$ be two rational numbers such that $p<r<s<p'$.

By the coupling Lemma~\ref{yal}, $\mathcal{G}^1_{p,q} \subset \mathcal{G}^1_{r,q}$ and $\mathcal{G}^2_{p',q} \subset \mathcal{G}^2_{s,q}$, whence
$$A\subset\bigcup_{
\begin{subarray}{c}
0\le r < s \le q \\
(r,s)\in\Q^2
\end{subarray}}
 \mathcal{G}^1_{r, q}\cap \mathcal{G}^2_{s, q}.$$
Then, it follows from the previous step that $A$ has probability $0$.

\vspace{0.2cm}
\underline{Step 3.} Proof of Theorem~\ref{THETHEtheorem}. Let $ q>0$.

\vspace{0.2cm}
Let $n\ge 1$, and let $E$ a finite subset of the set of real numbers $p\in[0,\min(q,\pcfleche))$ such that $\P(\mathcal{G}^1_{p, q}\cap \mathcal{G}^2_{p, q})\ge 1/n$.
By the previous step
$$\sum_{p\in E}\1_{\mathcal{G}^1_{p, q}\cap \mathcal{G}^2_{p, q}}\le 1,
\text{ which implies that }
\sum_{p\in E} \P(\mathcal{G}^1_{p, q}\cap \mathcal{G}^2_{p, q})\le 1.$$
Thus the set of $p$ such that $p\le  q$ and $\P(\mathcal{G}^1_{p, q}\cap \mathcal{G}^2_{p, q})\ge 1/n$ contains at most $n$ points, which proves the theorem.


\def\refname{References}
\bibliographystyle{plain}
\bibliography{ncd}

\begin{thebibliography}{10}

\bibitem{aizenman-barsky}
Michael Aizenman and David~J. Barsky.
\newblock Sharpness of the phase transition in percolation models.
\newblock {\em Comm. Math. Phys.}, 108(3):489--526, 1987.

\bibitem{CCGKS}
J.~T. Chayes, L.~Chayes, G.~R. Grimmett, H.~Kesten, and R.~H. Schonmann.
\newblock The correlation length for the high-density phase of {B}ernoulli
  percolation.
\newblock {\em Ann. Probab.}, 17(4):1277--1302, 1989.

\bibitem{DHB}
M.~Deijfen, O.~H{\"a}ggstr{\"o}m, and J.~Bagley.
\newblock A stochastic model for competing growth on {$\bold R\sp d$}.
\newblock {\em Markov Process. Related Fields}, 10(2):217--248, 2004.

\bibitem{deijfen}
Maria Deijfen.
\newblock Asymptotic shape in a continuum growth model.
\newblock {\em Adv. in Appl. Probab.}, 35(2):303--318, 2003.

\bibitem{DH}
Maria Deijfen and Olle H{\"a}ggstr{\"o}m.
\newblock Coexistence in a two-type continuum growth model.
\newblock {\em Adv. in Appl. Probab.}, 36(4):973--980, 2004.

\bibitem{DH-nonmon}
Maria Deijfen and Olle H{\"a}ggstr{\"o}m.
\newblock Nonmonotonic coexistence regions for the two-type richardson model on
  graphs.
\newblock {\em preprint}, 2005.

\bibitem{DS}
R.~Durrett and R.~H. Schonmann.
\newblock Large deviations for the contact process and two-dimensional
  percolation.
\newblock {\em Probab. Theory Related Fields}, 77(4):583--603, 1988.

\bibitem{GM-large}
Olivier Garet and R\'egine Marchand.
\newblock Large deviations for the chemical distance in supercritical bernoulli
  percolation.
\newblock {\em Preprint, available at
  \texttt{http://arxiv.org/abs/math.PR/0409317}}, 2004.

\bibitem{GM-coex}
Olivier Garet and R{\'e}gine Marchand.
\newblock Coexistence in two-type first-passage percolation models.
\newblock {\em Ann. Appl. Probab.}, 15(1A):298--330, 2005.

\bibitem{GM-fpppc}
Olivier Garet and Régine Marchand.
\newblock Asymptotic shape for the chemical distance and first-passage
  percolation on the infinite {B}ernoulli cluster.
\newblock {\em ESAIM Probab. Statist.}, 8:169--199, 2004.

\bibitem{Grimmett-Marstrand}
G.~R. Grimmett and J.~M. Marstrand.
\newblock The supercritical phase of percolation is well behaved.
\newblock {\em Proc. Roy. Soc. London Ser. A}, 430(1879):439--457, 1990.

\bibitem{Grimmett-book}
Geoffrey Grimmett.
\newblock {\em Percolation}, volume 321 of {\em Grundlehren der Mathematischen
  Wissenschaften [Fundamental Principles of Mathematical Sciences]}.
\newblock Springer-Verlag, Berlin, second edition, 1999.

\bibitem{grimmett-kesten}
Geoffrey Grimmett and Harry Kesten.
\newblock First-passage percolation, network flows and electrical resistances.
\newblock {\em Z. Wahrsch. Verw. Gebiete}, 66(3):335--366, 1984.

\bibitem{Haggstrom-Pemantle-1}
Olle H{\"a}ggstr{\"o}m and Robin Pemantle.
\newblock First passage percolation and a model for competing spatial growth.
\newblock {\em J. Appl. Probab.}, 35(3):683--692, 1998.

\bibitem{Haggstrom-Pemantle-2}
Olle H{\"a}ggstr{\"o}m and Robin Pemantle.
\newblock Absence of mutual unbounded growth for almost all parameter values in
  the two-type {R}ichardson model.
\newblock {\em Stochastic Process. Appl.}, 90(2):207--222, 2000.

\bibitem{hoffman}
Christopher Hoffman.
\newblock Coexistence for {R}ichardson type competing spatial growth models.
\newblock {\em Ann. Appl. Probab.}, 15(1B):739--747, 2005.

\bibitem{KL}
George Kordzakhia and Steven~P. Lalley.
\newblock A two-species competition model on $\mathcal{Z}^d$.
\newblock {\em Stochastic Process. Appl.}, 115(5):781--796, 2005.

\bibitem{marchand}
R.~Marchand.
\newblock Strict inequalities for the time constant in first passage
  percolation.
\newblock {\em Ann. Appl. Probab.}, 12(3):1001--1038, 2002.

\bibitem{neuhauser}
Claudia Neuhauser.
\newblock Ergodic theorems for the multitype contact process.
\newblock {\em Probab. Theory Related Fields}, 91(3-4):467--506, 1992.

\bibitem{vdb-kes}
J.~van~den Berg and H.~Kesten.
\newblock Inequalities for the time constant in first-passage percolation.
\newblock {\em Ann. Appl. Probab.}, 3(1):56--80, 1993.

\end{thebibliography}


\end{document}